\RequirePackage{etoolbox}
\csdef{input@path}{{style/}{graphics/}}
\documentclass[aos,MSNbibl,seceqn,secthm,dvips]{arximspdf}
\usepackage{multirow}
\usepackage{dcolumn}

\doi{10.1214/14-AOS1281} 
\volume{43}
\issue{1}
\pubyear{2015}
\firstpage{382}
\lastpage{421}
\docsubty{FLA}

\makeatletter
\newcommand{\ds}{\displaystyle}
\newcolumntype{d}[1]{D{.}{.}{#1}}
\newcommand{\Trr}{\operatorname{Tr}}

\newcommand{\rrvert}{\vert}

\newcommand{\llvert}{\vert}
\newtheorem{lem}[thm]{Lemma}
\newtheorem{cor}[thm]{Corollary}
\newproclaim{defi}[thm]{Definition}
\newproclaim{rem}[thm]{Remark}
\newproclaim{con}[thm]{Condition}
\newproclaim{Exa}{Example}[section]
\makeatother

\begin{document}
\begin{frontmatter}

\title{Universality for the largest eigenvalue of sample covariance
matrices with general population}
\runtitle{Edge universality for sample covariance matrices}

\begin{aug}
\author[A]{\fnms{Zhigang}~\snm{Bao}\corref{}\thanksref{T1}\ead
[label=e1]{zhigangbao@zju.edu.cn}},
\author[B]{\fnms{Guangming}~\snm{Pan}\thanksref{T2}\ead
[label=e2]{gmpan@ntu.edu.sg}} \and
\author[C]{\fnms{Wang}~\snm{Zhou}\ead
[label=e3]{stazw@nus.edu.sg}\thanksref{T3}\ead
[label=u1,url]{http://www.sta.nus.edu.sg/\textasciitilde stazw/}}
\runauthor{Z. Bao, G. Pan and W. Zhou}
\affiliation{Zhejiang University, Nanyang Technological University
and\\
National University~of~Singapore}
\address[A]{Z. Bao\\
Department of Mathematics\\
Zhejiang University\\
Hangzhou,  310027\\
P. R. China\\
\printead{e1}}
\address[B]{G. Pan\\
Division of Mathematical Sciences \\
School of Physical\\
\quad and Mathematical Sciences\\
Nanyang Technological University\\
Singapore 637371\\
Singapore\\
\printead{e2}}
\address[C]{W. Zhou\\
Department of Statistics\\
\quad and Applied Probability\\
National University of Singapore\\
Singapore 117546\\
Singapore\\
\printead{e3}\\
\printead{u1}}
\end{aug}
\thankstext{T1}{Supported in part by NSFC Grant
11071213, ZJNSF Grant R6090034 and SRFDP Grant 20100101110001.}
\thankstext{T2}{Supported in part by the Ministry of Education,
Singapore, under Grant  ARC 14/11.}
\thankstext{T3}{Supported in part by the Ministry of Education,
Singapore, under Grant  ARC 14/11, and by a~Grant
R-155-000-131-112 at the National University of Singapore.}

\received{\smonth{12} \syear{2013}}
\revised{\smonth{10} \syear{2014}}

\begin{abstract}
This paper is aimed at deriving the universality of the largest
eigenvalue of a class of high-dimensional real or complex sample
covariance matrices of the form $\mathcal{W}_N=\Sigma^{1/2} XX^*\Sigma
^{1/2}$. Here, $X=(x_{ij})_{M,N}$ is an $M\times N$ random matrix with
independent entries $x_{ij},1\leq i\leq M, 1\leq j\leq N$ such that
$\mathbb{E}x_{ij}=0$, $\mathbb{E}|x_{ij}|^2=1/N$. On dimensionality, we
assume that $M=M(N)$ and $N/M\rightarrow d\in(0,\infty)$ as
$N\rightarrow\infty$. For a class of general deterministic
positive-definite $M\times M$ matrices $\Sigma$, under some additional
assumptions on the distribution of $x_{ij}$'s, we show that the
limiting behavior of the largest eigenvalue of $\mathcal{W}_N$ is
universal, via pursuing a Green function comparison strategy raised
in [\textit{Probab. Theory Related Fields} \textbf{154} (2012) 341--407,
\textit{Adv. Math.} \textbf{229}  (2012) 1435--1515] by Erd\H{o}s, Yau and Yin for Wigner matrices
and extended by Pillai and Yin
[\textit{Ann. Appl. Probab.} \textbf{24} (2014) 935--1001]
to sample covariance
matrices in the null case ($\Sigma=I$). Consequently, in the standard
complex case ($\mathbb{E}x_{ij}^2=0$), combing this universality
property and the results known for Gaussian matrices obtained by El
Karoui in [\textit{Ann. Probab.} \textbf{35} (2007) 663--714]
(nonsingular case) and Onatski in
[\textit{Ann. Appl. Probab.} \textbf{18} (2008) 470--490]
(singular case), we show that after an appropriate
normalization the largest eigenvalue of $\mathcal{W}_N$ converges
weakly to the type 2 Tracy--Widom distribution $\mathrm{TW}_2$.
Moreover, in the real case, we show that when $\Sigma$ is spiked with a
fixed number of subcritical spikes, the type 1 Tracy--Widom limit
$\mathrm{TW}_1$ holds for the normalized largest eigenvalue of $\mathcal
{W}_N$, which extends a result of F\'{e}ral and P\'{e}ch\'{e} in
[\textit{J. Math. Phys.} \textbf{50} (2009) 073302]
to the scenario of nondiagonal $\Sigma$ and more generally
distributed $X$. In summary, we establish the Tracy--Widom type
universality for the largest eigenvalue of generally distributed sample
covariance matrices under quite light assumptions on $\Sigma$.
Applications of these limiting results to statistical signal detection
and structure recognition of separable covariance matrices are also discussed.
\end{abstract}

\end{frontmatter}
%

\section{Introduction}\label{sec1}

In recent decades, researchers working on
multivariate analysis have a growing interest in data with large size
arising from various fields such as genomics, image processing,
microarray, proteomics and finance, to name but a few. The classical
setting of \textit{fixed $p$ and large $n$} may lose its validity in
tackling some statistical problems for high-dimensional data, due to
the so-called \textit{curse of dimensionality}. As a feasible and
useful way in dealing with high-dimensional data, the spectral analysis
of high-dimensional sample covariance matrices has attracted
considerable interests among statisticians, probabilitists and
mathematicians. Study toward the eigenvalues of sample covariance
matrices traces back to the works of Fisher \cite{Fisher1939}, Hsu \cite
{Hsu1939} and Roy \cite{Roy1939}, and becomes flourishing after the
seminal work of Mar\v{c}enko and Pastur \cite{MP1967},
in which the authors established the limiting spectral distribution (MP
type distribution) for a class of sample covariance matrices, under the
setting that $p$ and $n$ are comparable. Since then, a lot of research
has been devoted to understanding the asymptotic properties of various
spectral statistics of high-dimensional sample covariance matrices. One
can refer to the monograph of Bai and Silverstein \cite{BS2009} for a
comprehensive summary and detailed references.

In this paper, we will focus on the limiting behavior of the largest
eigenvalue of a class of high-dimensional sample covariance matrices,
which is of great interest naturally from the principal component
analysis point of view. The largest eigenvalue has been commonly used
in hypothesis testing problems on the structure of the population
covariance matrix. Not trying to be comprehensive, one can refer to
\cite{BDMN2011,Johnstone2001,Onatski2011,PG2009,Karoui2007} for
instance. We also refer to the review paper of Johnstone \cite{Johnstone2007} for further reading on the statistical motivations of
the study on the largest eigenvalue of sample covariance matrices.
Precisely, we will consider the sample covariance matrix of the form
%
\begin{equation}
\label{88888} \mathcal{W}=\mathcal{W}_N:=\Sigma^{1/2}XX^*
\Sigma^{1/2}, \qquad X=(x_{ij})_{M,N},
\end{equation}
where $\{x_{ij}:=x_{ij}(N),1\leq i\leq M:=M(N),1\leq j\leq N\}$ is a
collection of independent real or complex variables such that
\[
\mathbb{E}x_{ij}=0, \qquad \mathbb{E}|x_{ij}|^2=N^{-1}.
\]
We call $\mathcal{W}_N$ a \textit{standard} complex sample covariance
matrix if there also exists
\[
\mathbb{E}x_{ij}^2=0,\qquad 1\leq i\leq M,1\leq j\leq N.
\]
In addition, $\Sigma:=\Sigma_N$ is assumed to be an $M\times M$
positive-definite matrix. In particular, if the columns of $X$ are
independently drawn from $\mathbf{h}/\sqrt{N}$ for some random vector
$\mathbf{h}$ possessing covariance matrix $I$, $\mathcal{W}$ can then
be viewed as the sample covariance matrix of $N$ observations of the
random vector $\Sigma^{1/2}\mathbf{h}$. Conventionally, we call
$\mathcal{W}$ a Wishart matrix if $x_{ij}$'s are Gaussian.
As is well known now, the limiting distributions of the largest
eigenvalues for classical high-dimensional random matrices were
originally discovered by Tracy and Widom in \cite{TW1994,TW1996} for
Gaussian Wigner ensembles G(O$/$U$/$S)E, thus named as the Tracy--Widom law
of type~$\beta$ ($\beta=1,2,4$ for GOE, GUE, GSE, resp.), denoted
by $\mathrm{TW}_{\beta}$ hereafter. The analogs in the context of sample
covariance matrices with $\Sigma=I$ were carried out by Johansson
\cite{Johansson2000} and Johnstone \cite{Johnstone2001}. More specifically,
the $\mathrm{TW}_2$ and $\mathrm{TW}_1$ limits were established for the
largest eigenvalues of standard complex and real null Wishart matrices
in \cite{Johansson2000} and \cite{Johnstone2001}, respectively.

For the nonnull population covariance matrix, that is, $\Sigma\neq I$,
much work has been devoted to the so-called spiked model, introduced by
Johnstone in \cite{Johnstone2001}. We say  $\mathcal{W}$ is spiked when
a few eigenvalues of $\Sigma$ are not equal to $1$. On the spiked
Wishart models, one can refer to \cite{BBP2006} for the standard
complex case and \cite{BV20111,BV2011,Mo2012,Paul2007,WD2012} for the
real case. However, in most cases, $\Sigma$ has more complicated
structures. In this paper, a more general setting on $\Sigma$ stated in
(iii) of Condition \ref{con.1.1} below will be employed, whereby
El Karoui showed in \cite{Karoui2007} that the $\mathrm{TW}_2$ limit
holds for the standard complex nonnull Wishart matrices when $d>1$
(nonsingular case), followed by Onatski's extension to the singular
case ($0<d\leq1$) in \cite{Onatski2008}.

With the above mentioned limiting results for the Wishart matrices at
hand, a conventional sequel in the Random Matrix Theory is to establish
the so-called universality property for generally distributed sample
covariance matrices, which states that the limiting behavior of an
eigenvalue statistic usually does not depend on the details of the
distribution of the matrix entries. The universality property of the
extreme eigenvalues is usually referred to as \textit{edge universality}.
Specifically, for sample covariance matrices in the null case, the
Tracy--Widom law has been established for $\mathcal{W}$ under very
general assumptions on the distribution of $X$. The readers can refer
to \cite{Soshnikov2002,Peche2009,WK2012,PY2012} for some
representative developments on this topic. For generally distributed
spiked models, the universality property was also partially obtained in
\cite{BY2008} and \cite{FP2009}. Especially, in the latter, the authors
proved that $\mathrm{TW}_1$ also holds for real spiked sample
covariance matrices with a finite number of \textit{subcritical} spikes
(see Corollary~\ref{thm.1205.1} for definition).

In this paper, armed with the condition on $\Sigma$, that is, Condition
\ref{con.1.1}(iii), we will prove the universality of the largest
eigenvalues of $\mathcal{W}$. It will be clear that such a class of
$\Sigma$ contains those spiked population covariance matrices with a
finite number of subcritical spikes, and goes far beyond. This work can
therefore be viewed as a substantial generalization of the Tracy--Widom
type edge universality, verified for the null case in \cite{PY2012} and
\cite{WK2012}, to a class of nonnull sample covariance matrices under
quite light assumptions on $\Sigma$. A direct consequence of such a
universality property, together with the results in \cite{Karoui2007}
and \cite{Onatski2008}, is that the $\mathrm{TW}_2$ also holds for
generally distributed standard complex $\mathcal{W}$ under our setting
on $\Sigma$; see Corollary~\ref{thm.1.1}. Moreover, by combining the
aforementioned result in \cite{FP2009}, we can also show that $\mathrm{TW}_1$ holds for real sample covariance matrices with spiked $\Sigma$
containing a fixed number of subcritical spikes; see Corollary~\ref{thm.1205.1}.
Note that $\Sigma$ is required to be diagonal in \cite{FP2009} and all
odd order moments of $x_{ij}$'s are assumed to vanish. We stress here,
our result can remove these restrictions. Both Corollary~\ref{thm.1.1}
and Corollary~\ref{thm.1205.1} can be used in high-dimensional
statistical inference then. In Section~\ref{sec2}, we will introduce two
applications, namely \textit{Presence of signals in the correlated
noise} and \textit{one-sided identity of separable covariance matrix}.
Related numerical simulations will also be conducted.

In the sequel, we will start by introducing some notation and then
present our main results. Subsequently, we will give a brief
introduction of the so-called \textit{Green function comparison
strategy}, and then sketch our new inputs for treating the general
setting of $\Sigma$.

\subsection{Main results}\label{sec11}
Henceforth, we will denote by $\lambda_{n}(A)\leq\cdots\leq\lambda
_2(A)\leq\lambda_1(A)$ the ordered eigenvalues of an $n\times n$
Hermitian matrix $A$.
For simplicity, we set the dimensional ratio
\[
d_N:=N/M\rightarrow d\in(0,\infty)\qquad \mbox{as } N\rightarrow
\infty.
\]
The empirical spectral distribution (ESD) of $\Sigma$ is
\[
H_N(\lambda):=\frac{1}M\sum_{i=1}^M
\mathbf{1}_{\{\lambda_i(\Sigma)\leq
\lambda\}}
\]
and that of $\mathcal{W}$ is
\[
\underline{F}_N(\lambda):=\frac{1}M\sum
_{i=1}^M\mathbf{1}_{\{\lambda
_i(\mathcal{W})\leq\lambda\}}.
\]
Here and throughout the following, $\mathbf{1}_{\mathbb{S}}$ represents
the indicator function of the event~$\mathbb{S}$. In addition, we will
need a crucial parameter $\mathbf{c}:=\mathbf{c}(\Sigma,N,M)\in
[0,1/\lambda_1(\Sigma))$ satisfying the equation
%
\begin{equation}
\label{0.0}
\int \biggl(\frac{\lambda\mathbf{c}}{1-\lambda\mathbf{c}} \biggr)^2\,dH_N(
\lambda)=d_N.
\end{equation}
It is elementary to check that the solution to (\ref{0.0}) in
$[0,1/\lambda_1(\Sigma))$ is unique.
With the above notation at hand, we can state our main condition as follows.

\begin{con} \label{con.1.1}
Throughout the paper, we need the following
conditions.
\begin{longlist}[(ii)]
\item[(i)] (On dimensionality). We assume that there are some positive
constants $c_1$~and~$C_1$ such that $c_1<d_N<C_1$.

\item[(ii)] (On $X$). We assume that $\{x_{ij}:=x_{ij}(N),1\leq i\leq M,1\leq
j\leq N\}$ is a collection of independent real or complex variables
such that
$\mathbb{E}x_{ij}=0$ and $\mathbb{E}|x_{ij}|^2=N^{-1}$.
Moreover, we assume that $\sqrt{N}x_{ij}$'s have a uniform
subexponential tail, that is, there exists some positive constant $\tau
_0$ independent of $i,j,N$ such that for sufficiently large $t$, one has
%
\begin{equation}\label{k.30}
\mathbb{P}\bigl(|\sqrt{N}x_{ij}|\geq t\bigr)\leq\tau_0^{-1}
\exp\bigl(-t^{\tau_0}\bigr).
\end{equation}

\item[(iii)] (On $\Sigma$). We assume that $\liminf_N\lambda_M(\Sigma)>0$,
$\limsup_N\lambda_1(\Sigma)<\infty$ and
%
\begin{equation}\label{itit}
\limsup_N\lambda_1(\Sigma) \mathbf{c}<1.
\end{equation}
\end{longlist}
\end{con}

Besides,\vspace*{-3pt} we also need the following ad hoc terminology.

\begin{defi}[(Matching to order $k$)] \label{def.1.2} Let $X^{\mathbf
{u}}=(x^{\mathbf{u}}_{ij})_{M,N}$ and $X^{\mathbf{v}}=(x^{\mathbf
{v}}_{ij})_{M,N}$ be\vspace*{1pt} two matrices satisfying (ii) of Condition \ref
{con.1.1}. We say\vspace*{1pt} $X^{\mathbf{u}}$ matches $X^{\mathbf{v}}$ to order
$k$, if for all $1\leq i\leq M,1\leq j\leq N$ and nonnegative integers
$l,m$ with $l+m\leq k$, there exists
%
\begin{eqnarray}
&& \mathbb{E}\bigl(\Re\bigl(\sqrt{N}x^{\mathbf{u}}_{ij}
\bigr)^l\Im\bigl(\sqrt{N}x^{\mathbf
{u}}_{ij}
\bigr)^m\bigr)
\nonumber
\\[-8pt]
\label{99999}
\\[-8pt]
\nonumber
&& \qquad=\mathbb{E}\bigl(\Re\bigl(\sqrt{N}x^{\mathbf{v}}_{ij}
\bigr)^l\Im\bigl(\sqrt {N}x^{\mathbf{v}}_{ij}
\bigr)^m\bigr)+O\bigl(e^{-(\log N)^C}\bigr)
\end{eqnarray}
with some positive constant $C>1$.
Alternatively, if (\ref{99999}) holds, we also say that $\mathcal
{W}^{\mathbf{u}}$ matches $\mathcal{W}^{\mathbf{v}}$ to order $k$,
where $\mathcal{W}^{\mathbf{u}}=\Sigma^{1/2}X^{\mathbf{u}}(X^{\mathbf
{u}})^*\Sigma^{1/2}$ and $\mathcal{W}^{\mathbf{v}}=\Sigma
^{1/2}X^{\mathbf{v}}(X^{\mathbf{v}})^*\Sigma^{1/2}$.
\end{defi}\vspace*{-3pt}

Our main theorem on edge universality of $\mathcal{W}$ can be
formulated as\vspace*{-3pt} follows.
%
\begin{thm}[(Universality for both real and complex cases)] \label
{thm.0.0} Suppose that two sample covariance matrices $\mathcal
{W}^{\mathbf{u}}=\Sigma^{1/2}X^{\mathbf{u}}(X^{\mathbf{u}})^*\Sigma
^{1/2}$ and
$\mathcal{W}^{\mathbf{v}}=\Sigma^{1/2}X^{\mathbf{v}}(X^{\mathbf
{v}})^*\Sigma^{1/2}$ satisfy Condition \ref{con.1.1}, where $X^{\mathbf
{u}}:=(x^{\mathbf{u}}_{ij})_{M,N}$ and $X^{\mathbf{v}}:=(x^{\mathbf
{v}}_{ij})_{M,N}$.\vspace*{-3pt} Let
%
\begin{equation}\label{1.2}
\lambda_r=\frac{1}{\mathbf{c}} \biggl(1+d_N^{-1}
\int\frac{\lambda\mathbf
{c}}{1-\lambda\mathbf{c}}\,dH_N(\lambda) \biggr).
\end{equation}
Then for sufficiently large $N$ and any real number $s$ which may
depend on $N$, there exist some positive constants $\varepsilon,\delta>0$
such that
%
\begin{eqnarray}
&& \mathbb{P}\bigl(N^{2/3}\bigl(\lambda_1\bigl(
\mathcal{W}^{\mathbf{u}}\bigr)-\lambda_r\bigr)  \leq
s-N^{-\varepsilon}\bigr)-N^{-\delta}\nonumber\\
&& \label{0.1.2.3.2.1}\qquad \leq\mathbb{P}\bigl(N^{2/3}
\bigl(\lambda_1\bigl(\mathcal{W}^{\mathbf{v}}\bigr)-\lambda
_r\bigr)\leq s\bigr)
\\
&&\qquad \leq  \mathbb{P}\bigl(N^{2/3}\bigl(\lambda _1\bigl(
\mathcal{W}^{\mathbf{u}}\bigr)-\lambda_r\bigr)\leq
s+N^{-\varepsilon}\bigr)+N^{-\delta}\nonumber
\end{eqnarray}
if one of the following two additional conditions holds:
\begin{longlist}[\textbf{A}:]
\item[\textbf{A}:] $\Sigma$ is diagonal and $\mathcal{W}^{\mathbf{u}}$ matches
$\mathcal{W}^{\mathbf{v}}$ to order $2$.

\item[\textbf{B}:] $\mathcal{W}^{\mathbf{u}}$ matches $\mathcal{W}^{\mathbf{v}}$
to order $4$.
\end{longlist}
\end{thm}

\begin{rem} \label{rem.20141005}
Theorem~\ref{thm.0.0} can be extended to the case of joint distribution
of the largest $k$ eigenvalues for any fixed positive integer $k$,
namely, for any real numbers $s_1,\ldots, s_k$ which may depend on $N$,
there exist some positive constants $\varepsilon,\delta>0$ such that
\begin{eqnarray*}
&& \mathbb{P}\bigl(N^{2/3}\bigl(\lambda_1\bigl(
\mathcal{W}^{\mathbf{u}}\bigr)-\lambda_r\bigr)\leq
s_1-N^{-\varepsilon},\ldots,\\
&&\quad N^{2/3}\bigl(
\lambda_k\bigl(\mathcal{W}^{\mathbf
{u}}\bigr)-\lambda_r
\bigr)\leq s_k-N^{-\varepsilon}\bigr)-N^{-\delta}
\\
&&\qquad\leq  \mathbb{P}\bigl(N^{2/3}\bigl(\lambda_1\bigl(
\mathcal{W}^{\mathbf{v}}\bigr)-\lambda _r\bigr)\leq
s_1,\ldots, N^{2/3}\bigl(\lambda_k\bigl(
\mathcal{W}^{\mathbf{v}}\bigr)-\lambda _r\bigr)\leq
s_k\bigr)
\nonumber
\\
&&\qquad\leq  \mathbb{P}\bigl(N^{2/3}\bigl(\lambda_1\bigl(
\mathcal{W}^{\mathbf{u}}\bigr)-\lambda _r\bigr)\leq
s_1+N^{-\varepsilon},\ldots,\\
&& \hspace*{45pt}\qquad N^{2/3}\bigl(
\lambda_k\bigl(\mathcal{W}^{\mathbf
{u}}\bigr)-\lambda_r
\bigr)\leq s_k+N^{-\varepsilon}\bigr)+N^{-\delta}.
\end{eqnarray*}
Such an extension can be realized through a parallel discussion as that
for the null case in \cite{PY2012}. One can refer to \cite{PY2012} for
more details. Here, we do not reproduce it.
\end{rem}

Combining Theorem~\ref{thm.0.0} with Theorem~1 of \cite{Karoui2007} and
Proposition~2 of \cite{Onatski2008} yields the following more concrete
result in the standard complex case ($\mathbb{E}x_{ij}^2=0$).

\begin{cor}[(Tracy--Widom limit for the standard complex case)]
\label{thm.1.1}
Let $\mathcal{W}_N^{\mathbf{g}_\mathbb{C}}$ be a standard
complex Wishart matrix and $\mathcal{W}_N$ be a general standard
complex sample covariance matrix. Assume that both of them satisfy
Condition \ref{con.1.1}. Denoting
%
\begin{equation}\label{201405050001}
\sigma^3=\frac{1}{\mathbf{c}^3} \biggl(1+d_N^{-1}
\int \biggl(\frac{\lambda
\mathbf{c}}{1-\lambda\mathbf{c}} \biggr)^3\,dH_N(\lambda)
\biggr),
\end{equation}
we have
\[
N^{2/3} \biggl(\frac{\lambda_1(\mathcal{W}_N)-\lambda_r}{\sigma} \biggr)\Longrightarrow
\mathrm{TW}_2
\]
if either $\Sigma$ is diagonal or $\mathcal{W}_N$ matches $\mathcal
{W}^{\mathbf{g}_\mathbb{C}}_N$ to order $4$.
\end{cor}

\begin{rem} \label{rem.43211234}
According to Remark~\ref{rem.20141005},
we also have the fact that the joint distribution of
\[
\biggl(\frac{\lambda_1(\mathcal{W}_N)-\lambda_r}{\sigma},\ldots, \frac
{\lambda_k(\mathcal{W}_N)-\lambda_r}{\sigma} \biggr)
\]
converges weakly to the $k$-dimensional joint $\mathrm{TW}_2$.
\end{rem}

For real sample covariance matrices, putting our Theorem~\ref{thm.0.0}
and Theorem~1.6 of \cite{FP2009} together, we can get the following corollary.

\begin{cor}[(Tracy--Widom limit for the real spiked case)]\label
{thm.1205.1} Suppose that $\mathcal{W}_N$ is a real sample covariance
matrix satisfying \textup{(i)} and \textup{(ii)} of Condition \ref{con.1.1}. Let $r$ be
some given positive integer. Assume that $\Sigma$ is spiked in the
sense that $\lambda_1(\Sigma)\geq\cdots\geq\lambda_r(\Sigma)\geq\lambda
_{r+1}(\Sigma)=\cdots=\lambda_M(\Sigma)=1$. Moreover, the $r$ spikes
$\lambda_{i}(\Sigma),i=1,\ldots,r$ are fixed (independent of $N$) and
subcritical, that is,
$
\lambda_1(\Sigma)<1+(\sqrt{d})^{-1}$.
Let $\mathcal{W}_N^{\mathbf{g}_{\mathbb{R}}}$ be a real Wishart matrix
with population covariance matrix $\Sigma$. Then in the scenario of
$d\in[1,\infty)$ (i.e., nonsingular case), we have
\[
N^{2/3} \biggl(\frac{\lambda_1(\mathcal{W}_N)-\lambda_r}{\sigma} \biggr)\Longrightarrow
\mathrm{TW}_1
\]
if either $\Sigma$ is diagonal or $\mathcal{W}_N$ matches $\mathcal
{W}^{\mathbf{g}_\mathbb{R}}_N$ to order $4$, where $\sigma$ is defined
in (\ref{201405050001}).
In addition, we have
%
\begin{equation}\label{1206.12}
\hspace*{6pt}\lambda_r=\bigl(1+d_N^{-1/2}
\bigr)^2+O\bigl(N^{-1}\bigr), \qquad\sigma =d_N^{-1/2}
\bigl(1+d_N^{1/2}\bigr)^{4/3}+o(1).
\end{equation}
\end{cor}

\begin{rem}\label{rem.1206.1} Analogously, under the assumption of
Theorem~\ref{thm.1205.1} we can get that the joint distribution of the
first $k$ normalized eigenvalues converges weakly to the
$k$-dimensional joint $\mathrm{TW}_1$.
\end{rem}

\begin{rem}
Lemma~\ref{lem.2014050501} below will show that the special spiked
$\Sigma$ with a fixed number of subcritical spikes satisfies (iii) of
Condition \ref{con.1.1}. It is known that if there is any spike on or
above the critical value $1+(\sqrt{d})^{-1}$, the limiting distribution
of the largest eigenvalue will not be the classical Tracy--Widom law
any more, assuming $r$ is fixed. One can refer to \cite{BBP2006} and
\cite{BV2011} for such a phase transition phenomenon for the standard
complex and real cases, respectively. Such a fact reflects that (iii)
of Condition \ref{con.1.1} is quite light for the Tracy--Widom type
universality to hold.
\end{rem}

\begin{rem} We conjecture that the $\mathrm{TW}_1$ law holds for all
$\Sigma$ satisfying (iii) of Condition \ref{con.1.1} and the
restriction on the nonsingular case is also artificial. However, as far
as we know, only \cite{FP2009} can provide us the reference matrix to
use the universality property in the real case. This is why we just
focus on the special real spiked sample covariance matrices here.
Nevertheless, these restrictions do not conceal the generality of the
universality result (Theorem~\ref{thm.0.0}) itself even in the real case.
\end{rem}

\subsection{Basic notions}\label{sec12}
We define the $N\times N$ matrix
\[
W=W_N:=X^*\Sigma X
\]
which shares the same nonzero eigenvalues with $\mathcal{W}$. Denoting
the ESD of $W_N$ by~$F_{N}$, we see
%
\begin{equation}\label{1208.6}
F_N=d_N^{-1}\underline{F}_N+
\bigl(1-d_N^{-1}\bigr)\mathbf{1}_{[0,\infty)}.
\end{equation}
If there is some deterministic distribution $H$ such that
$H_N\Longrightarrow H$ as $N\rightarrow\infty$,
it is well known that there are deterministic distributions $F_{d,H}$
and $\underline{F}_{d,H}$ such that
$
F_{N}\Longrightarrow F_{d,H}$ and $\underline{F}_{N}\Longrightarrow
\underline{F}_{d,H}
$
in probability. One can refer to \cite{Bai1999} or \cite{BS2009} for
detailed discussions. Analogous to (\ref{1208.6}), we have the relation
%
\begin{equation}\label{101010}
F_{d,H}=d^{-1}\underline{F}_{d,H}+
\bigl(1-d^{-1}\bigr)\mathbf{1}_{[0,\infty)}.
\end{equation}

For any distribution function $D$, its Stieltjes transform $m_D(z)$ is
defined by
\[
m_D(z)=\int\frac{1}{\lambda-z}\,dD(\lambda)
\]
for all $z\in\mathbb{C}^+:=\{\omega\in\mathbb{C},\Im\omega>0\}$. And
for any square matrix $A$, its Green function is defined by
$G_A(z)=(A-zI)^{-1}, z\in\mathbb{C}^+$. For convenience, we will
denote the Green functions of $W_N$ and $\mathcal{W}_N$, respectively, by
\[
G(z)=G_N(z):=(W_N-z)^{-1}\hspace*{10pt} \mbox{and}\hspace*{10pt}
\mathcal{G}(z)=\mathcal {G}_N(z):=(\mathcal{W}_N-z)^{-1},
\hspace*{18pt} z\in\mathbb{C}^+.
\]
The Stieltjes transforms of $F_N$ and $\underline{F}_N$ will be denoted
by $m_N(z)$ and $\underline{m}_N(z)$, respectively. By definitions,
obviously one has
\[
m_N(z)=\frac{1}N \Trr G(z),\qquad \underline{m}_N(z)=
\frac{1}M \Trr \mathcal{G}(z).
\]
Here, we draw attention to the basic relation $\Trr G(z)-\Trr \mathcal{G}(z)=(M-N)/z$.
Actually, what really pertains to our discussion in the sequel is the
nonasymptotic version of $F_{d,H}$ which can
be obtained via replacing $d$ and $H$ by $d_N$ and $H_N$ in $F_{d,H}$,
and thus will be denoted by $F_{d_N,H_N}$. More precisely, $F_{d_N,H_N}$
is the corresponding distribution function of the Stieltjes transform
$m_{d_N,H_N}(z):=m_{F_{d_N,H_N}}(z)\in\mathbb{C}^{+}$ satisfying the
following self-consistent equation:
%
\begin{equation}
\label{1.3}
\hspace*{10pt}\quad m_{d_N,H_N}(z)=\frac{1}{-z+d_N^{-1}\int {t}/({tm_{d_N,H_N}(z)+1})\,dH_N(t)},\qquad z\in\mathbb{C}^{+}.
\end{equation}

Analogously, we can define the nonasymptotic versions of $\underline
{F}_{d,H}$ and its Stieltjes transform, denoted by $\underline
{F}_{d_N,H_N}$ and $\underline{m}_{d_N,H_N}(z)$, respectively. Then the
$N$-dependent version of (\ref{101010}) is
%
\begin{equation}\label{1208.5}
F_{d_N,H_N}=d_N^{-1}\underline{F}_{d_N,H_N}+
\bigl(1-d_N^{-1}\bigr)\mathbf{1}_{[0,\infty)}.
\end{equation}

For simplicity, we will briefly use the notation
\begin{eqnarray*}
m_0(z)&:=& m_{d_N,H_N}(z),\qquad \underline{m}_0(z):=
\underline {m}_{d_N,H_N}(z),\\
F_0 &:= & F_{d_N,H_N},\qquad
\underline {F}_0:=\underline{F}_{d_N,H_N}
\end{eqnarray*}
in the sequel.

It has been discussed in \cite{SC1995} by Silverstein and Choi that
$F_0$ has a continuous derivative $\rho_0$ on $\mathbb{R}\setminus\{0\}
$ and the rightmost boundary of the support of $\rho_0$ is $\lambda_r$
defined in (\ref{1.2}), that is, $\lambda_r=\inf\{x\in\mathbb
{R}\dvtx F_0(x)=1\}$.
Moreover, the parameter $\mathbf{c}$ defined by (\ref{0.0}) satisfies
$\mathbf{c}=-\lim_{z\in\mathbb{C}^{+}\rightarrow\lambda_r}m_0(z)$.

\subsection{Sketch of the proof route}\label{sec13}
As mentioned above, Theorem~\ref{thm.0.0} can be viewed as a substantial generalization of the edge
universality for the null sample covariance matrices provided in \cite
{PY2012}. However, the general machinery in \cite{PY2012}, with the
so-called Green function comparison approach at the core, still works
well even for general nonnull case. The Green function comparison
strategy was raised in the series of work
\cite{EYY2011,EYY2012,EYY20122} on the local eigenvalue statistics of
Wigner matrices originally, and has shown its strong applicability on
some other random matrix models or statistics; see \cite{EKYY2012,BPZ2012,PY20121} for its variants for sample covariance and
correlation matrices and see \cite{TV2011} for an application on random
determinant. We also refer to the survey \cite{Erdos2011} for an overview.

To be specific, the preliminary heuristic of the Green function
comparison strategy for our objective can be roughly explained as
follows. At first, the distribution function of $\lambda_1(\mathcal
{W})$ can actually be approximated from above and from below by the
expectations of two functionals of the Stieltjes transform $m_N(z)$,
that is, the normalized trace of the Green function $G_N(z)$; see (\ref
{2014050501}) below. Hence, the comparison between the distributions of
the largest eigenvalues of $\mathcal{W}^{\mathbf{u}}$ and $\mathcal
{W}^{\mathbf{v}}$ can then be reduced to the comparison between the
expectations of the functionals of the Green functions. For the latter,
a replacement method inherited from the classical Lindeberg swapping
process (see \cite{Lindeberg1922}) can be employed. Together with the
expansion formula of the Green function, such a replacement method can
effectively lead to the universality property.

A main technical tool escorting the Green function comparison process
is the so-called \textit{strong law of local eigenvalue density}, which
asserts that the limiting spectral law is even valid on short intervals
which contain only $N^{\varepsilon}$ eigenvalues for any constant
$\varepsilon>0$. Such a limiting law on microscopic scales was
developed in a series of work \cite{ESY2009,ESY20092,ESY2010,EYY20122} for Wigner matrices originally and was shown to be crucial in
recent work on universality problems of local eigenvalue statistics,
one can refer to \cite{ESY2011,EYY2012,TV2011} for instance. For our
purpose, we will need a \textit{strong local MP type law around} $\lambda_r$, which was established in our recent paper \cite{BPZ2013} and is
recorded as Theorem~\ref{thm.3.90} below. The companion work \cite{BPZ2013} initiates the project of edge universality and provides
essential technical inputs for the Green function comparison process.
However, the strong law of local eigenvalue density is also of interest
in its own right.

To lighten the notation, we make the convention $E=\Re z$ and $\eta=\Im
z$ hereafter. And we also denote
\[
\Delta(z):=\Sigma^{1/2}\mathcal{G}(z)\Sigma^{1/2}
\]
for simplicity. It will be seen that, in our comparison process, we
need to control the magnitude of the entries of $\Delta(z)$ in the
regime $|E-\lambda_r|\leq N^{-2/3+\varepsilon}$ and $\eta
=N^{-2/3-\varepsilon}$ for some small positive constant $\varepsilon$.
This issue turns out to be a new difficulty due to the complexity of
$\Sigma$. We handle this main technical task for diagonal and
nondiagonal $\Sigma$ via substantially different approaches, which are
sketched as follows.

Clearly, when $\Sigma$ is diagonal, we can turn to bound the entries of
$\mathcal{G}(z)$ instead. Invoking the spectral decomposition [see the
first inequality of (\ref{201404271}) below, e.g.], the desired
bound can be obtained via providing (1): an accurate description of the
locations of the eigenvalues; (2): an upper bound for the eigenvector
coefficients. It will be clear that (1) can be transformed into the
strong local MP type law which has already been established. Toward (2),
we will prove the so-called \textit{delocalization property}, which
states the eigenvector coefficients are of order $O(N^{-1/2+\varepsilon
})$ typically. The delocalization property was first derived in
\cite{ESY2009} and improved in the series of papers
\cite{ESY20092,ESY2010,EYY20122} for Wigner matrices, and extended to sample covariance
matrices in the null case in \cite{TV2012,WK2012,ESYY2012,PY2012}.
Here, we extend the delocalization property to $\mathcal{W}$ for those
eigenvectors corresponding to the eigenvalues around $\lambda_r$.

However, for the nondiagonal $\Sigma$, we need to focus on the entries
of $\Delta(z)$ themselves. Fortunately, it turns out that only the
diagonal entries $\Delta_{kk}(z)$ should be bounded if we are
additionally granted in the comparison process that two ensembles match
to order 4. To this end, we can start from the spectral decomposition
again [see (\ref{2014050311}) below, e.g.]. Analogous to the
diagonal case, we could provide (1$'$): an accurate description of the
locations of the eigenvalues; (2$'$): an upper bound for $(\Sigma
^{1/2}\mathbf{u}_i\mathbf{u}_i^*\Sigma^{1/2})_{kk}$. Observe that (1$'$)
is just the same as (1) for diagonal $\Sigma$, actually can also be
ensured by the strong local MP type law. However, (2$'$) requires some
totally novel ideas. More details in Section~\ref{sec5} will show that the
spectral decomposition equality (\ref{2014050311}) can also be applied
in a converse direction, to wit, with a bound on $\Delta_{kk}(z_0)$ for
some appropriately chosen $z_0:=E_0+\mathbf{i}\eta_0$, one can
actually obtain a bound for $(\Sigma^{1/2}\mathbf{u}_i\mathbf
{u}_i^*\Sigma^{1/2})_{kk}$ in turn. For $\eta_0=N^{-2/3+\varepsilon}\gg
\eta$, we will perform a novel bounding scheme for $\Delta_{kk}(z_0)$,
based on the Schur complement and the concentration inequalities on
quadratic forms (Lemma~\ref{lem.k.40}). Then, by the bound on $\Delta
_{kk}(z_0)$ one can get a bound on $(\Sigma^{1/2}\mathbf{u}_i\mathbf
{u}_i^*\Sigma^{1/2})_{kk}$, which together with (1$'$) implies the
desired bound on $\Delta_{kk}(z)$. The choice of $\eta_0\gg N^{-2/3}$
will be technically necessary for our bounding scheme on $\Delta
_{kk}(z_0)$. Therefore, we adopt such a roundabout way to bound $\Delta
_{kk}(z)$, owing to the fact that $\eta\ll N^{-2/3}$ is unaffordable
for a direct application of our bounding scheme based on the Schur
complement and the concentration inequalities.

\subsection{Notation and organization}\label{sec14}
Throughout the paper, we use the notation $O(\cdot)$ and $o(\cdot)$ in
the conventional sense. As usual, $C,C_1,C_2$ and $C'$ stand for some
generic positive constants whose values may differ from line to line.
We say $x \sim y$
if there exist some positive constants $C_1$ and $C_2$ such that
$C_1|y|\leq|x|\leq C_2|y|$. Generally, for two functions
$f(z),g(z)\dvtx
\mathbb{C}\rightarrow\mathbb{C}$, we say $f(z)\sim g(z)$
if there exist some positive constants $C_1$ and $C_2$ independent of
$z$ such that $C_1|g(z)|\leq|f(z)|\leq C_2|g(z)|$. Moreover,
$\Vert A\Vert_{\mathrm{op}}$ and $\Vert A\Vert_{\mathrm{HS}}$ represent the operator norm and
Hilbert--Schmidt norm of a matrix $A$, respectively, and $\Vert\mathbf{u}\Vert$ is the $L_2$ norm of a vector $\mathbf{u}$. We use $\mathbf{i}$
to denote the imaginary unit to release $i$ which will be frequently
used as index or subscript. In addition, we conventionally denote by
$\mathbf{e}_i$ the vector with all 0's except for a 1 in the $i$th
coordinate and by $\mathbf{1}$ the vector with 1 in each coordinate.
The dimensions of these vectors are usually obvious according to the
context thus just omitted from the notation. $\mathbf{0}_{\alpha\times
\beta}$ will be used to represent the $\alpha\times\beta$ null matrix
which will be abbreviated to $\mathbf{0}_{\alpha}$ if $\alpha=\beta$.
In addition, we adopt the notation in \cite{PY2012} to set the
frequently used parameter
\[
\varphi:=\varphi_N=(\log N)^{\log\log N}.
\]
For $\zeta>0$, we say that an event $\mathbb{S}$ holds with $\zeta
$-high probability if there is some positive constant $C$ such that for
sufficiently large $N$,
\[
\mathbb{P}(\mathbb{S})\geq1-N^{C}\exp\bigl(-\varphi^{\zeta}
\bigr).
\]

We conclude this section by stating its organization. In Section~\ref{sec2}, we
will introduce some applications of our main results in
high-dimensional statistical inference, and some related simulations
will be conducted. Then we will turn to the theoretical part.
In Section~\ref{sec3}, we will recall the properties of $m_0(z)$ and the strong
local MP type law around $\lambda_r$ established in \cite{BPZ2013} as
the preliminaries of our proofs for the main results. In Section~\ref{sec4}, we
will use the strong local MP type law and a Green function comparison
approach to prove Theorem~\ref{thm.0.0}, Corollaries \ref{thm.1.1} and~\ref{thm.1205.1}. Section~\ref{sec5} will be devoted to the aforementioned
argument of bounding the entries of $\Delta(z)$.

\section{Applications and simulations}\label{sec2}

In this section, we introduce some applications of our universality
results in high-dimensional statistical inference, and conduct related
simulations to check the quality of the approximations of our limiting
laws and discuss their utility in the concrete hypothesis testing
problems. We remark here, though
Corollary~\ref{thm.1205.1} and Remark~\ref{rem.1206.1} are only stated
for the case of $d\geq1$, we will also perform the simulations for the
case of $d<1$.

\subsection{Applications}\label{sec21}
\begin{itemize}
\item \textit{Presence of signals in the correlated noise.}
\end{itemize}

Consider an $M$-dimensional signal-plus-noise vector $\mathbf{y}:=A\mathbf{s}+\Sigma_a^{1/2}\mathbf{z}$ and its $N$ i.i.d. samples, namely
\[
\mathbf{y}_i=A\mathbf{s}_i+\Sigma_a^{1/2}
\mathbf{z}_i,\qquad i=1,\ldots,N,
\]
where $\mathbf{s}$ is a $k$-dimensional real or complex mean zero
signal vector with covariance matrix $S$; $\mathbf{z}$ is an
$M$-dimensional real or complex random vector with independent mean
zero and variance one coordinates; $A$ is an $M\times k$ deterministic
matrix which is of full column rank and $\Sigma_a$ is an $M\times M$
deterministic positive-definite matrix. We call $\Sigma_a^{1/2}\mathbf
{z}$ the noise vector. Moreover, the signal vector and the noise vector
are assumed to be independent. Set the matrices $Z_N=[\mathbf
{z}_1,\ldots,\mathbf{z}_N]$ and $Y_N=[\mathbf{y}_1,\ldots,\mathbf
{y}_N]$. Denoting the covariance matrix of $\mathbf{y}$ by $R$, we can
get by assumption that
\[
R=\mathit{ASA}^T+\Sigma_a.
\]
Such a model stems from several statistical signal processing problems,
and is used commonly in various fields such as wireless communications,
bioinformatics and machine learning, to name a few. We refer to Kay
\cite{Kay1998} for a comprehensive overview. A fundamental target is to
detect signals via data. Thus the very first step is to know whether
there is any signal present, that is, $k=0$ versus $k\geq1$. Once
signals are detected, one can take a step further to estimate the
number $k$. Under the high-dimensional setting, Nadakuditi and Edelman
in \cite{NE2008}, and Bianchi et al. in \cite{BDMN2011} considered
respectively to detect signals in the \textit{white} Gaussian noise,
that is, $\Sigma_a=I$ (or more generally, $\Sigma_a=cI$ with some
positive number $c$) and $\mathbf{z}$ is Gaussian. Also under the
Gaussian assumption on the noise, Nadakuditi and Silverstein in \cite{NS2010} considered this detection problem when the noise may be
correlated, that is, $\Sigma_a$ may not be a multiple of $I$. We also
refer to the very recent work of Vinogradova, etc. \cite{VCH2014} for
the case of correlated noise.
Our aim is to test, for generally distributed and correlated noise
$\Sigma_a^{1/2}\mathbf{z}$, whether there is no signal present. Thus
our hypothesis testing problem can be stated as
\[
(\mathbf{Q}_a)\mbox{:}\qquad\mathbf{H}_0\mbox{:}\qquad
k=0 \quad \mbox{vs.} \quad \mathbf{H}_1\mbox{:}\qquad
k\geq1.
\]
\begin{itemize}
\item \textit{One-sided identity of separable covariance matrices.}
\end{itemize}

Consider the data model of the form
\[
\mathcal{Y}_N=\Sigma_b^{1/2}
\mathcal{Z}_NT^{1/2},
\]
where $\mathcal{Z}_N$ is an $M\times N$ random matrix sharing the same
distribution as $Z_N$ in the previous problem, $T$ is an $N\times N$
deterministic positive-definite matrix and $\Sigma_b$ is an $M\times M$
deterministic positive-definite matrix. $N^{-1}\mathcal{Y}_N\mathcal
{Y}_N^*$ is then called the \textit{separable covariance matrix} which
is widely used for handling the spatiotemporal sampling data. Such a
nomenclature is owing to the fact that the vectorization of the data
matrix $\mathcal{Y}_N$ has a separable covariance $\Sigma_b\otimes T$.
The spectral properties of $N^{-1}\mathcal{Y}_N\mathcal{Y}_N^*$ have
been widely investigated in some recent work under the high-dimensional
setting, for example, one can refer to \cite{PS2009,Karoui2009,Zhang2006,WP2013}. Without loss of generality, we regard $T$ as the
temporal covariance matrix and $\Sigma_b$ as the spatial covariance
matrix. In this paper, we are interested in whether the temporal
identity (i.e., $T=I$) holds. Formally, we are concerned with the
following hypothesis testing:
\[
(\mathbf{Q}_b)\mbox{:}\qquad \mathbf{H}_0\mbox{:}
\qquad T=I \quad \mbox{vs.}\quad\mathbf{H}_1\mbox{:}\qquad
T\neq I.
\]
Actually,\hspace*{1.5pt} we can consider to test whether $T=T_0$ for any given
positive-definite~$T_0$, since considering the renormalized data matrix
$\mathcal{Y}_NT_0^{-1/2}$ we can recover the testing problem $\mathbf{Q}_b$.
A similar testing problem with $T$ replaced by $\Sigma_b$ can also be
considered. We call this kind of hypothesis testing problem \textit{one-sided identity test} for the separable covariance matrix.
\begin{itemize}
\item \textit{Onatski's statistics.}
\end{itemize}

Note that under $\mathbf{H}_0$ of either $\mathbf{Q}_a$ or $\mathbf
{Q}_b$, the involved sample covariance matrix $N^{-1}Y_NY_N^*$ or
$N^{-1}\mathcal{Y}_N\mathcal{Y}_N^*$ is of the form $\mathcal{W}$
defined in (\ref{88888}).
It is then natural to construct our test statistics for $\mathbf{Q}_a$
and $\mathbf{Q}_b$ from the largest eigenvalues of $N^{-1}Y_NY_N^*$ and
$N^{-1}\mathcal{Y}_N\mathcal{Y}_N^*$, respectively, such that our
universality results can be employed under $\mathbf{H}_0$. For
simplicity, we will use $\mathcal{W}$ to represent either
$N^{-1}Y_NY_N^*$ or $N^{-1}\mathcal{Y}_N\mathcal{Y}_N^*$ under $\mathbf
{H}_0$, that is, we will regard $(\Sigma, X)$ as $(\Sigma_a, Z_N/\sqrt
{N})$ and $(\Sigma_b,\mathcal{Z}_N/\sqrt{N})$ when we refer to $\mathbf
{Q}_a$ and $\mathbf{Q}_b$, respectively.

At first glance, it is natural to choose the normalized largest
eigenvalue as our test statistic. Unfortunately, in the real system,
$\Sigma$ is usually\vadjust{\goodbreak} unknown. Hence, a~general result like Corollary~\ref
{thm.1.1}, where the parameters $\lambda_r$ and $\sigma$ depend on
$\Sigma$, cannot be used directly if no information of $\Sigma$ is
known priori. To eliminate the unknown parameters $\lambda_r$ and
$\sigma$, we adopt the strategy used by Onatski in \cite{Onatski2009,Onatski20092}. More specifically, we will use the statistics
\[
\mathbf{T}_a=\frac{\lambda_1(Y_NY_N^*)-\lambda_2(Y_NY_N^*)}{\lambda
_2(Y_NY_N^*)-\lambda_3(Y_NY_N^*)} \quad\mbox{and}\quad \mathbf{T}_b=
\frac{\lambda_1(\mathcal{Y}_N\mathcal{Y}_N^*)-\lambda
_2(\mathcal{Y}_N\mathcal{Y}_N^*)}{\lambda_2(\mathcal{Y}_N\mathcal
{Y}_N^*)-\lambda_3(\mathcal{Y}_N\mathcal{Y}_N^*)}
\]
for $\mathbf{Q}_a$ and $\mathbf{Q}_b$, respectively. In the sequel, we
will call $\mathbf{T}_a$ and $\mathbf{T}_b$ \textit{Onatski's
statistics}. Note that under $\mathbf{H}_0$, $\mathbf{T}_a$ and
$\mathbf{T}_b$ possess the same limiting distribution, determined by
the joint $\mathrm{TW}_\beta$ laws, mentioned in Remarks \ref
{rem.43211234} and \ref{rem.1206.1}. An obvious advantage of $\mathbf
{T}_a$ or $\mathbf{T}_b$ is that its limiting distribution is
independent of $\lambda_r$ and $\sigma$ under $\mathbf{H}_0$, which
makes it asymptotically pivotal. Moreover, though the explicit formula
for the limiting distribution function of Onatski's statistic under
$\mathbf{H}_0$ is unavailable currently, one can approximate it via
simulation, by generating the eigenvalues from high-dimensional GOE
(resp., GUE) in the real (resp., complex) case. We will describe such
an approximation in detail in the subsequent simulation study.

\subsection{Simulations}\label{sec22}
\begin{itemize}
\item \textit{Accuracy of approximations for TW laws.}
\end{itemize}
We conduct some numerical simulations to check the accuracy of the
distributional approximations in Corollaries \ref{thm.1.1} and \ref
{thm.1205.1}, under various settings of $(M,N)$, $\Sigma$ and the
distribution of $X$. Firstly, for each pair of $(M,N)$, we generate an
observation from $M\times M$ Haar distributed random orthogonal matrix
and denote it by $U:=U(M,N)$. To get such a $U$, we can generate in
Matlab an $M\times M$ Gaussian matrix $\mathbf{G}$ with i.i.d. $N(0,1)$
entries, and let $U=\mathbf{G}(\mathbf{G}^*\mathbf{G})^{-1/2}$ which is
well defined with probability 1; refer to Section~7.1 of \cite{Eaton1989} for instance. Then we will fix this $U$ for each pair of
$(M,N)$ as a deterministic orthogonal matrix. Next, we set some
scenarios of $\Sigma$ in Corollaries \ref{thm.1.1} and \ref
{thm.1205.1}. To this end, we define
\[
D_c:=\operatorname{diag}(\underbrace{1,\ldots,1}_{\lfloor M/2\rfloor
},
\underbrace{2,\ldots,2}_{M-\lfloor M/2\rfloor}),\qquad D_r:=\operatorname{diag}
\biggl(1+\frac{(\sqrt{d_N})^{-1}}{2},1,\ldots,1\biggr)
\]
and choose $\Sigma$ to be some similar forms of $D_c$ and $D_r$ in
Corollaries \ref{thm.1.1} and \ref{thm.1205.1}, respectively. More
specifically, we will use the following four choices of population
covariance matrix $\Sigma$, denoted by
\begin{eqnarray*}
\Sigma(c,1) &:=& D_c,\qquad \Sigma(c,2):=UD_cU^*,\\
\Sigma (r,1) &:=& D_r, \qquad\Sigma(r,2):=UD_rU^*.
\end{eqnarray*}
Here, $U$ is the orthogonal matrix generated priori.

Now, we state our choice for the distribution of $X$. For simplicity,
we set $h_{ij}:=\sqrt{N}x_{ij}$ for $i=1,\ldots,M$ and $j=1,\ldots,N$,
and choose all these $(h_{ij})$'s to be i.i.d. For standard complex
Gaussian case, the numerical performance of the limiting law in
Corollary~\ref{thm.1.1} has been assessed; see Tables~1 and 2 of \cite{Karoui2007}. Here, we use a discrete distribution in our simulation
study. Specifically, let $\mathfrak{s}_1$ and $\mathfrak{s}_2$ be two
i.i.d. variables with the distribution
\[
\mathfrak{u}=\tfrac{1}{12}\bolds{\delta}_{-2}+
\tfrac
{4}{25}\bolds{\delta}_{-1}+\tfrac{13}{24}
\bolds{\delta}_0+\tfrac{16}{75}\bolds{
\delta}_{{3}/{2}}+\tfrac{1}{600}\bolds{\delta}_4,
\]
where $\bolds{\delta}_a$ represents the Dirac measure at $a$. It
is elementary to check that the first four moments of $\mathfrak{u}$
are the same as those of $N(0,1)$. Now we choose $h_{11}$ for the
standard complex and real cases respectively as
\[
h_{11}(c)\stackrel{d}=\tfrac{1}{\sqrt{2}}(\mathfrak{s}_1+
\mathbf {i}\mathfrak{s}_2) \quad\mbox{and}\quad h_{11}(r)\stackrel{d}=\mathfrak{s}_1
\]
and denote the corresponding $X$ by $X(c)$ and $X(r)$, respectively. We
conduct the simulations for the combinations $(\Sigma(c,1),X(c))$,
$(\Sigma(c,2),X(c))$, $(\Sigma(r,1),X(r))$ and $(\Sigma(r,2),X(r))$
under various settings of $(M,N)$. The results are provided in Table~\ref{t1}.
It can be seen, in each case, the approximation is satisfactory even
for relatively small $M$ and $N$.

\begin{table}
\tabcolsep=0pt
\tablewidth=\textwidth
\caption{Simulated quantiles for four pairs of $(\Sigma,X)$. The cases
of $(\Sigma(r,1),X(r))$, $(\Sigma(r,2), X(r))$, $(\Sigma(c,1), X(c))$
and $(\Sigma(c,2), X(2))$ are titled by {R1}, {R2}, {C1} and {C2},
respectively, for\vspace*{-2pt} simplicity}\label{t1}
{\fontsize{8.3}{10.3}{\selectfont
\begin{tabular*}{\tablewidth}{@{\extracolsep{\fill}}ld{2.4}d{1.4}d{1.4}d{1.4}d{1.4}d{1.4}d{1.4}d{1.4}d{1.4}@{}}
\hline
\multicolumn{1}{@{}l}{$\bolds{(\Sigma,X)}$}& \multicolumn{1}{c}{\textbf{Percentile}}& \multicolumn{1}{c}{${\mathbf{TW}}_{\mathbf{1}}$}&
\multicolumn{1}{c}{$\bolds{30\times30}$}&
\multicolumn{1}{c}{$\bolds{60\times 60}$}&
\multicolumn{1}{c}{$\bolds{100\times 100}$}&
\multicolumn{1}{c}{$\bolds{80\times 20}$}&
\multicolumn{1}{c}{$\bolds{20 \times 80}$}&
\multicolumn{1}{c}{$\bolds{100 \times 400}$} &
\multicolumn{1}{c@{}}{$\bolds{2}\bolds{\ast} \mathbf{SE}$}\\
\hline
R1 &
-3.9000& 0.0100& 0.0053& 0.0087& 0.0114& 0.0075& 0.0076& 0.0115&
0.0020\\
&-3.1800& 0.0500& 0.0479& 0.0523& 0.0566& 0.0493& 0.0580& 0.0601&
0.0040\\
&-2.7800& 0.1000& 0.1070& 0.1151& 0.1151& 0.1099& 0.1197& 0.1192&
0.0060\\
&-1.9100& 0.3000& 0.3520& 0.3524& 0.3393& 0.3535& 0.3539& 0.3352&
0.0090\\
&-1.2700& 0.5000& 0.5762& 0.5674& 0.5457& 0.5725& 0.5714& 0.5388&
0.1000\\
&-0.5900& 0.7000& 0.7752& 0.7547& 0.7388& 0.7685& 0.7713& 0.7372&
0.0090\\
& 0.4500& 0.9000& 0.9345& 0.9260& 0.9214& 0.9347& 0.9359& 0.9171&
0.0060\\
& 0.9800& 0.9500& 0.9689& 0.9650& 0.9620& 0.9706& 0.9708& 0.9620&
0.0040\\
& 2.0200& 0.9900& 0.9943& 0.9929& 0.9931& 0.9938& 0.9952& 0.9905&
0.0020\\[3pt]
R2&
-3.9000& 0.0100& 0.0054& 0.0088& 0.0098& 0.0063& 0.0086& 0.0100&
0.0020\\
&-3.1800& 0.0500& 0.0497& 0.0533& 0.0552& 0.0504& 0.0556& 0.0572&
0.0040\\
&-2.7800& 0.1000& 0.1077& 0.1144& 0.1143& 0.1085& 0.1226& 0.1116&
0.0060\\
&-1.9100& 0.3000& 0.3617& 0.3460& 0.3363& 0.3470& 0.3707& 0.3392&
0.0090\\
&-1.2700& 0.5000& 0.5784& 0.5582& 0.5506& 0.5700& 0.5834& 0.5433&
0.1000\\
&-0.5900& 0.7000& 0.7714& 0.7568& 0.7503& 0.7731& 0.7765& 0.7404&
0.0090\\
& 0.4500& 0.9000& 0.9301& 0.9258& 0.9248& 0.9334& 0.9349& 0.9166&
0.0060\\
& 0.9800& 0.9500& 0.9658& 0.9649& 0.9630& 0.9671& 0.9704& 0.9605&
0.0040\\
& 2.0200& 0.9900& 0.9924& 0.9929& 0.9928& 0.9934& 0.9941& 0.9937&
0.0020\\[6pt]
\multicolumn{1}{@{}l}{$\bolds{(\Sigma,X)}$}&
\multicolumn{1}{c}{\textbf{Percentile}}&
\multicolumn{1}{c}{$\mathbf{TW}_{\mathbf{2}}$}&
\multicolumn{1}{c}{$\bolds{30 \times 30}$}&
\multicolumn{1}{c}{$\bolds{60 \times 60}$}&
\multicolumn{1}{c}{$\bolds{100 \times 100}$} &
\multicolumn{1}{c}{$\bolds{80 \times 20}$}&
\multicolumn{1}{c}{$\bolds{20 \times 80}$}&
\multicolumn{1}{c}{$\bolds{100 \times 400}$} &
\multicolumn{1}{c@{}}{$\bolds{2} \bolds{\ast}\mathbf{SE}$}\\
\hline
C1&
-3.7300& 0.0100& 0.0031& 0.0053& 0.0066& 0.0037& 0.0042& 0.0082&
0.0020\\
&-3.2000& 0.0500& 0.0266& 0.0377& 0.0363& 0.0319& 0.0326& 0.0396&
0.0040\\
&-2.9000& 0.1000& 0.0674& 0.0812& 0.0827& 0.0749& 0.0745& 0.0870&
0.0060\\
&-2.2700& 0.3000& 0.2573& 0.2728& 0.2819& 0.2772& 0.2648& 0.2819&
0.0090\\
&-1.8100& 0.5000& 0.4695& 0.4804& 0.4866& 0.4838& 0.4818& 0.4861& 0.1000\\
&-1.3300& 0.7000& 0.6913& 0.6963& 0.6950& 0.6942& 0.6928& 0.6936&
0.0090\\
&-0.6000& 0.9000& 0.9053& 0.9004& 0.9006& 0.9012& 0.9021& 0.9025&
0.0060\\
&-0.2300& 0.9500& 0.9549& 0.9506& 0.9489& 0.9521& 0.9525& 0.9531&
0.0040\\
& 0.4800& 0.9900& 0.9913& 0.9900& 0.9886& 0.9880& 0.9924& 0.9912&
0.0020\\[3pt]
C2&
-3.7300& 0.0100& 0.0021& 0.0056& 0.0066& 0.0032& 0.0050& 0.0071&
0.0020\\
&-3.2000& 0.0500& 0.0234& 0.0321& 0.0399& 0.0326& 0.0321& 0.0445&
0.0040\\
&-2.9000& 0.1000& 0.0642& 0.0754& 0.0852& 0.0781& 0.0746& 0.0926&
0.0060\\
&-2.2700& 0.3000& 0.2639& 0.2641& 0.2805& 0.2721& 0.2734& 0.2955&
0.0090\\
&-1.8100& 0.5000& 0.4745& 0.4756& 0.4858& 0.4874& 0.4864& 0.4933&
0.1000\\
&-1.3300& 0.7000& 0.6875& 0.6877& 0.6930& 0.7006& 0.6923& 0.6954&
0.0090\\
&-0.6000& 0.9000& 0.9008& 0.9012& 0.8988& 0.9055& 0.8994& 0.9028&
0.0060\\
&-0.2300& 0.9500& 0.9490& 0.9512& 0.9493& 0.9547& 0.9517& 0.9529&
0.0040\\
& 0.4800& 0.9900& 0.9899& 0.9905& 0.9894& 0.9917& 0.9891& 0.9903&
0.0020\\
\hline
\end{tabular*}
}}\vspace*{-6pt}
\tabnotetext[{}]{t1}{The simulation was done in Matlab. In each of the above four cases, we
generated 10{,}000 matrix $X$ with the distribution defined above, and
then calculated the largest eigenvalue of $\mathcal{W}$ and
renormalized it with the parameters $\lambda_r$ and $\sigma$ according
to Corollaries \ref{thm.1.1}~and~\ref{thm.1205.1}. In the column titled
``Percentile,'' we listed the quantiles of $\mathrm{TW}_\beta$ law for $\beta
=1,2$. Simulating 10{,}000 times gave us an empirical distribution of the
renormalized largest eigenvalue. And we stated the values of this
empirical distribution at the quantiles of the TW laws for various
pairs of $(M,N)=(30,30), (60,60), (100,100), (80, 20), (20, 80),
(100,400)$. The last column states the approximate standard errors
based on binomial sampling.}
\end{table}

Next, a natural question is, to what extent can we weaken the
assumptions imposed on $X$. Very recently, a necessary and sufficient
condition for the Tracy--Widom limit of Wigner matrix with i.i.d.
off-diagonal entries (up to symmetry) was established by Lee and Yin in
\cite{LY2014}, where the matrix entry is only required to have mean 0
and variance 1, and satisfies a tail condition slightly weaker than the
existence of the 4th moment. It is reasonable to conjecture a similar
moment condition is sufficient for the validity of Tracy--Widom laws
for sample covariance matrices. To give a numerical evidence for such a
conjecture, we also conduct some simulation for the largest eigenvalue
of $\mathcal{W}$ whose entries possess a symmetric Pareto distribution.
For simplicity, we only state the simulation result for the complex
case with $\Sigma=\Sigma(c,2)$. We choose $\mathfrak{p}_1$ and
$\mathfrak{p}_2$ to be i.i.d. variables with the symmetric Pareto
distribution whose density is given by $f(x)=\frac{9}{10}\sqrt{\frac
{3}{5}}|x|^{-6}$ when $|x|> \sqrt{\frac{3}{5}}$ and $0$ otherwise. It
is then elementary to see that
%
\begin{equation} \label{201404281}
h_{11}\stackrel{d}=\tfrac{1}{\sqrt{2}}(\mathfrak{p}_1+
\mathbf{i}\mathfrak {p}_2)
\end{equation}
has mean 0 and variance 1. Moreover, we see $\mathbb{E}|h_{11}|^4<\infty$. We denote by $X(P)$ the corresponding $X$. The simulation results
are stated in Table~\ref{t2}. It can be seen that the approximation is also
very good even for small $M$ and $N$.
%
\begin{table}
\tablewidth=\textwidth
\tabcolsep=0pt
\caption{Simulated quantiles for the case of $(\Sigma, X)=(\Sigma(c,2),
X(P))$ (CP for short)}\label{t2}
\begin{tabular*}{\tablewidth}{@{\extracolsep{\fill}}lccccccccc@{}}
\hline
$\bolds{(\Sigma, X)}$& \textbf{Percentile}& $\mathbf{TW}_{\mathbf{2}}$& $\bolds{30 \times 30}$&
$\bolds{60 \times 60}$&
$\bolds{100 \times  100}$& $\bolds{80 \times 20}$& $\bolds{20 \times 80}$&
$\bolds{100 \times 400}$ & $\bolds{2} \bolds{\ast}\mathbf{SE}$\\
\hline
CP&
$-3.7300$& 0.0100& 0.0016& 0.0043& 0.0062& 0.0035& 0.0044& 0.0088&
0.0020\\
&$-3.2000$& 0.0500& 0.0280& 0.0409& 0.0460& 0.0345& 0.0369& 0.0512&
0.0040\\
&$-2.9000$& 0.1000& 0.0776& 0.0987& 0.1037& 0.0862& 0.0894& 0.1069&
0.0060\\
&$-2.2700$& 0.3000& 0.3113& 0.3311& 0.3320& 0.3201& 0.3275& 0.3235&
0.0090\\
&$-1.8100$& 0.5000& 0.5517& 0.5476& 0.5507& 0.5603& 0.5628& 0.5347&
0.1000\\
&$-1.3300$& 0.7000& 0.7675& 0.7472& 0.7501& 0.7736& 0.7689& 0.7335&
0.0090\\
&$-0.6000$& 0.9000& 0.9392& 0.9303& 0.9228& 0.9364& 0.9380& 0.9176&
0.0060\\
&$-0.2300$& 0.9500& 0.9716& 0.9658& 0.9659& 0.9714& 0.9708& 0.9580&
0.0040\\
& \phantom{$-$}0.4800& 0.9900& 0.9932& 0.9909& 0.9921& 0.9928& 0.9928& 0.9912&
0.0020\\
\hline
\end{tabular*}
\tabnotetext[]{t1}{The simulation was taken analogously. We generated 10{,}000 matrix $X$
with $h_{11}$ following the distribution defined in (\ref{201404281}).
Each column has the same meaning as that in Table~\ref{t1}.}
\end{table}
\begin{itemize}
\item \textit{Size and power study for} $\mathbf{T}_a$ \textit{and} $\mathbf{T}_b$.
\end{itemize}
Now, we evaluate the sizes and powers of the statistics $\mathbf{T}_a$
and $\mathbf{T}_b$ for $\mathbf{Q}_a$ and $\mathbf{Q}_b$ respectively.
For simplicity, we only report the results for the real case here. Note
that, in the real case, we do not establish\vadjust{\goodbreak} the $\mathrm{TW}_1$ law for
general $\Sigma$ satisfying Condition \ref{con.1.1}(iii). However, in
the sequel, we will also perform the simulation for $\Sigma$ which is
not spiked, such as $\Sigma=\Sigma(c,1)$.
More specifically, we will focus on two settings
\[
\mbox{(I):}\qquad  \Sigma_a=\Sigma_b=\Sigma(r,1), \qquad Z_N\stackrel
{d}=\mathcal{Z}_N\stackrel{d}=\sqrt{N}X(r),
\]
and
\[
\mbox{(II):}\qquad
 \Sigma_a=\Sigma_b=\Sigma(c,1),\qquad  Z_N\stackrel
{d}=\mathcal{Z}_N\stackrel{d}=\sqrt{N}X(r).
\]
For $\mathbf{Q}_a$, we choose the alternative with some positive number
$\rho_a$ as
\[
\mathbf{H}_1(a,\rho_a)\mbox{:} \qquad k=1,\qquad A=\mathbf{e}'_1
\quad\mbox{and}\quad \mathbf{s}\sim N(0,\rho_a),
\]
where $\mathbf{e}_1$ is $M$-dimensional by the assumption on $A$.
For $\mathbf{Q}_b$, we choose two alternatives parameterized by $\rho
_b$ as
\[
\mathbf{H}_1(b,\rho_b,1)\mbox{:}\qquad T=I+\rho_b
\mathbf{e}_1\mathbf {e}'_1 \quad\mbox{and}\quad
\mathbf{H}_1(b,\rho_b,2)\mbox{:}\qquad T=I+\rho_b
\frac{1}{N}\mathbf{1}\mathbf{1}',
\]
where $\mathbf{e}_1$ and $\mathbf{1}$ are both $N$-dimensional by the
assumption on $T$. Under the setting (I),
for $\mathbf{H}_1(a,\rho_a)$, we set $\rho_{a}:=\tau(\sqrt{d_N})^{-1}$,
while for both $\mathbf{H}_1(b,\rho_b,1)$ and $\mathbf{H}_1(b,\rho
_b,2)$, we set $\rho_{b}:=\tau\sqrt{d_N}$ with some strength parameter
$\tau>0$. Under the setting (II), for $\mathbf{H}_1(a,\rho_a)$, we set
$\rho_{a}:=2\tau(\sqrt{d_N})^{-1}$, while for both $\mathbf{H}_1(b,\rho
_b,1)$ and $\mathbf{H}_1(b,\rho_b,2)$, we set $\rho_{b}:=2\tau\sqrt
{d_N}$ with some strength parameter $\tau>0$. We will choose $\tau
=0.5,4,6$ for each alternative above.

Now assuming that $\xi_1$, $\xi_2$ and $\xi_3$ have the joint $\mathrm{TW}_1$ distribution, we approximate the percentiles of the
distribution of $(\xi_1-\xi_2)/(\xi_2-\xi_3)$ as follows. We can
simulate 30{,}000 independent matrices from GOE of dimension 1000 and
numerically compute the ratio of the differences between the first and
the second and the second and the third eigenvalues for each matrix,
then we can get the percentiles of the empirical distribution of these
30,000 ratios. By doing the above in Matlab, we got that the
approximate 95th percentile of the distribution of $(\xi_1-\xi_2)/(\xi
_2-\xi_3)$ is 7.16. The nominal significant level of our tests is 5\%.
The results for the sizes are reported in Table~\ref{t3}, and the results for
the powers are reported in Table~\ref{t4} for setting (I) and Table~\ref{t5} for
(II), respectively. The small $\tau=0.5$ is tailored for corroborating
the following phenomenon, that is, \textit{an additive or multiplicative
finite rank perturbation may not cause significant change of the
largest eigenvalue of a sample covariance matrix when the strength of
the perturbation is weak enough.} This phenomenon has been explicitly
verified for the spiked sample covariance matrices, see the
aforementioned references on the spiked models \cite{BBP2006} and \cite{FP2009}. Our Corollary~\ref{thm.1205.1} also confirms it again.
However, for more complicated models such as $N^{-1}Y_NY_N^*$ and
$N^{-1}\mathcal{Y}_N\mathcal{Y}_N^*$ in our $\mathbf{Q}_a$ and $\mathbf
{Q}_b$, given general $\Sigma_a$ and $\Sigma_b$, the theoretical
discussions on this phenomenon with respect to various $A\mathbf{s}$
and $T$ are still open.
Under our choices of $\Sigma$, from the simulations we can see that
when $\tau=0.5$, the powers of both tests in various scenarios are very
poor. However, when $\tau$ is relatively large, our tests are reliable.
It can be seen from Tables~\ref{t4} and \ref{t5}, in the cases of $\tau=4$ or $6$,
the powers are satisfactory, especially when $N$ and $M$ are relatively large.
%
\begin{table}[t]
\tablewidth=\textwidth
\tabcolsep=0pt
\caption{Simulated sizes for settings \textup{(I)} and \textup{(II)}}\label{t3}
\begin{tabular*}{\tablewidth}{@{\extracolsep{\fill}}lccccccccc@{}}
\hline
\textbf{Setting} & $\bolds{30 \times 30}$& $\bolds{60 \times 60}$& $\bolds{100 \times 100}$&
$\bolds{80 \times 20}$&
$\bolds{80 \times 40}$& $\bolds{20 \times 80}$& $\bolds{40 \times 80}$& $\bolds{100 \times 400}$ & $\bolds{400 \times
200}$\\
\hline
\phantom{I}(I) &  0.0522& 0.0476& 0.0490& 0.0604& 0.0526& 0.0474& 0.0521& 0.0512&
0.0486\\
(II)& 0.0544& 0.0511& 0.0488& 0.0543& 0.0521& 0.0493& 0.0526& 0.0478&
0.0446\\
\hline
\end{tabular*}
\end{table}

\begin{table}[b]
\tablewidth=\textwidth
\tabcolsep=0pt
\caption{Simulated powers for $\mathbf{T}_a$ and $\mathbf{T}_b$ under
setting \textup{(I)}, $\tau$ is 0.5, 4 or 6}\label{t4}
{\fontsize{8.3}{10.3}{\selectfont
\begin{tabular*}{\tablewidth}{@{\extracolsep{\fill}}lcccccccccc@{}}
\hline
$\bolds{\tau}$& $\mathbf{H}_{\bolds{1}}$& $\bolds{30 \times 30}$& $\bolds{60 \times 60}$&
$\bolds{100 \times 100}$&
$\bolds{80 \times 20}$& $\bolds{80 \times 40}$&
$\bolds{20 \times 80}$&
$\bolds{40 \times 80}$&
$\bolds{100 \times 400}$ &
$\bolds{400 \times 200}$\\
\hline
0.5&
$\mathbf{H}_1(a,\rho_a)$& 0.0630& 0.0577& 0.0622& 0.0604& 0.0588&
0.0589& 0.0614& 0.0541& 0.0545\\
&$\mathbf{H}_1(b,\rho_b,1)$& 0.0541& 0.0516& 0.0497& 0.0540& 0.0521&
0.0533& 0.0522& 0.0488& 0.0482\\
&$\mathbf{H}_1(b,\rho_b,2)$& 0.0551& 0.0463& 0.0508& 0.0540& 0.0518&
0.0530& 0.0545& 0.0506& 0.0498\\[3pt]
4&
$\mathbf{H}_1(a,\rho_a)$& 0.4825& 0.6857& 0.8421& 0.5090& 0.6454&
0.5263& 0.6680& 0.9684& 0.9929\\
&$\mathbf{H}_1(b,\rho_b,1)$& 0.3932& 0.5775& 0.7507& 0.4243& 0.5475&
0.4262& 0.5488& 0.8998& 0.9816\\
&$\mathbf{H}_1(b,\rho_b,2)$& 0.3983& 0.5776& 0.7511& 0.4216& 0.5529&
0.4352& 0.5427& 0.8970& 0.9812\\[3pt]
6 &
$\mathbf{H}_1(a,\rho_a)$& 0.7319& 0.9089& 0.9807& 0.7434& 0.8830&
0.7861& 0.8980& 0.9999& 1.0000\\
&$\mathbf{H}_1(b,\rho_b,1)$& 0.6556& 0.8653& 0.9628& 0.7205& 0.8539&
0.6822& 0.8247& 0.9954& 0.9998\\
&$\mathbf{H}_1(b,\rho_b,2)$& 0.6623& 0.8647& 0.9608& 0.7235& 0.8477&
0.6888& 0.8277& 0.9955& 1.0000\\
\hline
\end{tabular*} }}
\end{table}

\begin{table}
\tablewidth=\textwidth
\tabcolsep=0pt
\caption{Simulated powers for $\mathbf{T}_a$ and $\mathbf{T}_b$ under
setting \textup{(II)}, $\tau$ is 0.5, 4 or 6}\label{t5}
{\fontsize{8.3}{10.3}{\selectfont\begin{tabular*}{\tablewidth}{@{\extracolsep{\fill}}lcccccccccc@{}}
\hline
$\bolds{\tau}$& $\mathbf{H}_{\bolds{1}}$& $\bolds{30 \times 30}$& $\bolds{60 \times 60}$&
$\bolds{100 \times 100}$&
$\bolds{80 \times 20}$& $\bolds{80 \times 40}$& $\bolds{20 \times 80}$&
$\bolds{40 \times 80}$& $\bolds{100 \times 400}$
& $\bolds{400 \times 200}$\\
\hline
0.5 &
$\mathbf{H}_1(a,\rho_a)$& 0.0583& 0.0588& 0.0534& 0.0649& 0.0557&
0.0568& 0.0554& 0.0536& 0.0512\\
&$\mathbf{H}_1(b,\rho_b,1)$& 0.0627& 0.0583& 0.0587& 0.0614& 0.0577&
0.0607& 0.0589& 0.0572& 0.0545\\
&$\mathbf{H}_1(b,\rho_b,2)$& 0.0620& 0.0593& 0.0532& 0.0614& 0.0618&
0.0590& 0.0565& 0.0513& 0.0502\\[3pt]
4 &
$\mathbf{H}_1(a,\rho_a)$& 0.4352& 0.6207& 0.7891& 0.5269& 0.6354&
0.3476& 0.5360& 0.8278& 0.9913\\
&$\mathbf{H}_1(b,\rho_b,1)$& 0.7650& 0.9258& 0.9870& 0.8454& 0.9277&
0.7237& 0.8776& 0.9974& 1.0000\\
&$\mathbf{H}_1(b,\rho_b,2)$& 0.7558& 0.9249& 0.9856& 0.8480& 0.9259&
0.7263& 0.8780& 0.9978& 1.0000\\[3pt]
6&
$\mathbf{H}_1(a,\rho_a)$& 0.9255& 0.9914& 0.9996& 0.9754& 0.9954&
0.9031& 0.9817& 1.0000& 1.0000\\
&$\mathbf{H}_1(b,\rho_b,1)$& 0.9594& 0.9978& 1.0000& 0.9854& 0.9984&
0.9595& 0.9951& 1.0000& 1.0000\\
&$\mathbf{H}_1(b,\rho_b,2)$& 0.9285& 0.9923& 0.9997& 0.9758& 0.9956&
0.9072& 0.9831& 1.0000& 1.0000\\
\hline
\end{tabular*}}}
\end{table}

\section{Square root behavior and local MP type law}\label{sec3}
In this section, we will record several main results proved in our
recent paper \cite{BPZ2013} which will serve as fundamental inputs for
the Green function comparison process. The main result established in
\cite{BPZ2013} is the aforementioned strong local MP type law around
$\lambda_r$; see Theorem~\ref{thm.3.90} below. As a
necessary input to the proof of the strong local MP type law around
$\lambda_r$ in \cite{BPZ2013}, the so-called square root behavior of
$m_0(z)$ has been derived therein, see Theorem~\ref{lem.2.2} below.
Then, as a direct consequence of the strong local MP type law around
$\lambda_r$, a nearly optimal convergence rate of $F_N(x)$ around
$\lambda_r$ has also been obtained, see Theorem~\ref{thm.1206.7} below.
All these results will play roles in our Green function comparison process.

In \cite{BPZ2013}, it has been shown that $\mathbf{c},\lambda_r\sim1$.
More precisely, there exist two positive constants $C_l\leq C_r$ such
that $\lambda_r\in[C_l/2,2C_r]$. Here, $C_l$ and $C_r$ can be chosen
appropriately such that $\lambda_1(\mathcal{W})\in[C_l,C_r]$ with $\zeta$-high probability for any given constant $\zeta>0$. We will always
write $z:=E+\mathbf{i}\eta$, and use the notation
\[
\kappa:=\kappa(z)=|E-\lambda_r|.
\]
We introduce for $\zeta\geq0$ two sets of $z$,
\begin{eqnarray*}
S(\zeta) &:=& \bigl\{z\in\mathbb{C}\dvtx  C_l\leq E\leq C_r,
\varphi^{\zeta
}N^{-1}\leq\eta\leq1\bigr\},
\\
S_r(\tilde{c},\zeta) &:=& \bigl\{z\in\mathbb{C}\dvtx
\lambda_r-\tilde{c}\leq E\leq C_r,
\varphi^{\zeta}N^{-1}\leq\eta\leq1\bigr\},
\end{eqnarray*}
where $\tilde{c}$ is some positive constant.

The first main result we need is a collection of some crucial
properties of $m_0(z)$, which are essentially guaranteed by (iii) of
Condition \ref{con.1.1}, and can be inferred from the square root
behavior of the limiting spectral density $\rho_0$ on its right edge
$\lambda_r$. Informally, we can call it \textit{square root behavior
of} $m_0(z)$.
%
\begin{thm}[(Square root behavior of $m_0(z)$, Lemma~2.3 of \cite{BPZ2013})] \label{lem.2.2}
Under Condition \ref{con.1.1}, there exists
some small but fixed positive constant $\tilde{c}$ such that the
following three statements hold.
\begin{longlist}[(iii)]
\item[(i)] For $z\in S(0)$, we have
\[
\bigl|m_0(z)\bigr| \sim 1;
\]

\item[(ii)]  For $z\in S_r(\tilde{c},0)$, we have
\[
\Im m_0(z)\sim
\cases{
\ds\frac{\eta}{\sqrt{\kappa+\eta}},&\quad \mbox{if }$E\geq\lambda_{r}+\eta$,\vspace*{6pt}
\cr
\ds\sqrt{\kappa+\eta},& \quad\mbox{if }$E\in[\lambda_r-\tilde{c},\lambda
_r+\eta)$;}
\]
\item[(iii)] For $z\in S_r(\tilde{c},0)$, we have
\[
\bigl|1+tm_0(z)\bigr|\geq\hat{c}\bigl(1+\lambda_{1}(
\Sigma)m_0(\lambda_r)\bigr)\geq c_0,\qquad
\forall t\in\bigl[\lambda_M(\Sigma),\lambda_1(\Sigma)
\bigr]
\]
for some small positive constants $\hat{c},c_0$ depending on $\tilde{c}$.
\end{longlist}
\end{thm}

The second necessary input is the strong local MP type law around
$\lambda_r$. To state it, we also need to recall some additional
notation from \cite{BPZ2013}. We denote by $\mathbf{x}_i$ the $i$th
column of $X$ and set
$\mathbf{r}_i=\Sigma^{1/2}\mathbf{x}_i$. We introduce the notation
$X^{(\mathbb{T})}$ to represent the $M\times(N-|\mathbb{T}|)$ minor of
$X$ obtained by deleting $\mathbf{x}_i$ from $X$ if $i\in\mathbb{T}$.
For convenience, $(\{i\})$ will be abbreviated to $(i)$. Denoting
\[
W^{(\mathbb{T})}=X^{(\mathbb{T})*}\Sigma X^{(\mathbb{T})},\qquad \mathcal
{W}^{(\mathbb{T})}=\Sigma^{1/2}X^{(\mathbb{T})}X^{(\mathbb{T})*}
\Sigma^{1/2},
\]
we can further set
\begin{eqnarray*}
G^{(\mathbb{T})}(z)&=& \bigl(W^{(\mathbb{T})}-z\bigr)^{-1},\qquad \mathcal
{G}^{(\mathbb{T})}(z)=\bigl(\mathcal{W}^{(\mathbb{T})}-z\bigr)^{-1},\\
m_N^{(\mathbb{T})}(z) &=& \frac{\Trr G^{(\mathbb{T})}(z)}{N},\qquad \underline
{m}_N^{(\mathbb{T})}(z)=\frac{\Trr \mathcal{G}^{(\mathbb{T})}(z)}{M}.
\end{eqnarray*}
We emphasize here, in the sequel, the names of indices of $X$ for
$X^{(\mathbb{T})}$ will be kept, that is, $X^{(\mathbb{T})}_{ij}=\mathbf
{1}_{\{j\notin\mathbb{T}\}} X_{ij}$. Correspondingly, we will denote
the $(i,j)$th entry of $G^{(\mathbb{T})}(z)$ by $G^{(\mathbb
{T})}_{ij}(z)$ for all $i,j\notin\mathbb{T}$. In addition, in light of
the discussion in \cite{BPZ2013} [see the truncation issue above (3.3)
therein], henceforth we can and do additionally assume\vspace*{-3pt} that
%
\begin{equation}\label{1206.5}
\max_{i,j}|\sqrt{N}x_{ij}|\leq(\log
N)^{C}
\end{equation}
with some sufficiently large positive constant $C$.
Then we have the following theorem.
%
\begin{thm}[(Strong local MP type law around $\lambda_r$, Theorem~3.2 of
\cite{BPZ2013})] \label{thm.3.90}
Let $\tilde{c}$ be the constant in
Theorem~\ref{lem.2.2}. Under Condition \ref{con.1.1} and assumption
(\ref{1206.5}), for any $\zeta>0$ there exists some constant $C_\zeta$ such\vspace*{-3pt} that
\begin{longlist}[(ii)]
\item[(i)]
%
\begin{equation}\label{9.1.1}
\bigcap_{z\in S_r(\tilde{c},5C_\zeta)} \biggl\{\bigl|m_N(z)-m_0(z)\bigr|
\leq\varphi ^{C_\zeta}\frac{1}{N\eta} \biggr\}
\end{equation}
holds with $\zeta$-high probability, and

\item[(ii)]
%
\begin{eqnarray}\label{9.1.2}
&&\bigcap_{z\in S_r(\tilde{c},5C_\zeta)} \biggl\{\max_{i\neq
j}\bigl|G_{ij}(z)\bigr|+
\max_{i}\bigl|G_{ii}(z)-m_0(z)\bigr|
\nonumber
\\[-8pt]
\\[-8pt]
\nonumber
&&\qquad\hspace*{70pt}\leq
\varphi^{C_\zeta} \biggl(\sqrt {\frac{\Im m_0(z)}{N\eta}}+\frac{1}{N\eta}
\biggr) \biggr\}
\end{eqnarray}
holds with $\zeta$-high probability.
\end{longlist}
\end{thm}

For our purpose, the following result concerning the convergence rate
of ESD around $\lambda_r$ is also needed, which can be essentially
derived from Theorem~\ref{thm.3.90}.
%
\begin{thm}[(Convergence rate around $\lambda_r$, Theorem~4.1 of \cite{BPZ2013})]\label{thm.1206.7}
Under Condition \ref{con.1.1} and the
assumption (\ref{1206.5}), for any $\zeta>0$ there exists a constant
$C_\zeta$ such that the following events hold with $\zeta$-high probability.
\begin{longlist}[(ii):]
\item[(i)]  For the largest eigenvalue $\lambda_1(\mathcal{W})$, there exists
\[
\bigl|\lambda_1(\mathcal{W})-\lambda_r\bigr|\leq
N^{-2/3}\varphi^{C_\zeta}.
\]

\item[(ii)] For any
$E_1,E_2\in[\lambda_r-\tilde{c},C_r]$,
there exists
%
\begin{equation}\label{5.1.900}
\bigl|\bigl(F_N(E_1)-F_N(E_2)
\bigr)-\bigl(F_0(E_1)-F_0(E_2)
\bigr)\bigr|\leq N^{-1}\varphi^{C_\zeta}.
\end{equation}
\end{longlist}
\end{thm}

In addition, we record the following concentration lemma on quadratic
forms, whose proof can be found in Appendix~B of \cite{EYY2012} for instance.

\begin{lem} \label{lem.k.40}Let $\mathbf{x}_i,\mathbf{x}_j,i\neq j$ be
two columns of the matrix $X$ satisfying \textup{(ii)} of Condition \ref
{con.1.1}. Then for any $M$-dimensional vector $\mathbf{b}$ and
$M\times M$ matrix $\mathbf{C}$ independent of $\mathbf{x}_i$ and
$\mathbf{x}_j$, the following three inequalities hold with $\zeta$-high
probability
\begin{longlist}[(ii)]
\item[(i)]
\[
 \biggl|\mathbf{x}_i^*\mathbf{C}\mathbf{x}_i-\frac{1}N\Trr\mathbf
{C}\biggr|\leq\frac{\varphi^{\tau\zeta}}{N}\Vert\mathbf{C}\Vert_{\mathrm{HS}},
\]
\item[(ii)]
\[
\bigl|\mathbf{x}_i^*\mathbf{C}\mathbf{x}_j\bigr|\leq\frac{\varphi^{\tau\zeta
}}{N}\Vert\mathbf{C}\Vert_{\mathrm{HS}},
\]
\item[(iii)]
\[
\bigl|\mathbf{b}^*\mathbf{x}_i\bigr|\leq\frac{\varphi
^{\tau\zeta}}{\sqrt{N}}\Vert\mathbf{b}\Vert.
\]
Here, $\tau:=\tau(\tau_0)> 1$ is some positive constant [see \textup{(ii)} of
Condition \ref{con.1.1} for $\tau_0]$. Let $\mathfrak{X}_i$ be the
conjugate transpose of the $i$th row of the matrix $X$ for $i=1,\ldots,
M$. If we replace $\mathbf{x}_i,\mathbf{x}_j$ by $\mathfrak{X}_i$ and
$\mathfrak{X}_j$ respectively, the above three inequalities also hold
if $\mathbf{b}$ is an $N$-dimensional vector and $\mathbf{C}$ is an
$N\times N$ matrix which are both independent of $\mathfrak{X}_i$ and
$\mathfrak{X}_j$.
\end{longlist}
\end{lem}

Finally, regarding the $\Vert\cdot\Vert_{\mathrm{HS}}$ norm of a Green function, we
will frequently need the fact that for any Hermitian matrix $A$, there is
%
\begin{eqnarray}
\bigl\Vert(A-z)^{-1}\bigr\Vert_{\mathrm{HS}}^2 &=& \Trr |A-z|^{-2}
= \Trr
(A-z)^{-1}(A-\bar{z})^{-1}
\nonumber
\\[-8pt]
\label{144201}\\[-8pt]
\nonumber
&=& \frac
{1}{\eta}\Im
\Trr(A-z)^{-1},
\end{eqnarray}
which can be verified easily by the spectral decomposition.
\section{Universality for the largest eigenvalue}\label{sec4}
With some bounds on the entries of $\Sigma^{1/2}\mathcal{G}\Sigma
^{1/2}$ granted (see Lemma~\ref{lem.11.11.105} below), we can
successfully prove our main results in this section via pursuing a
Green function comparison strategy tailored for edge universality,
which is analogous to those in \cite{PY2012,PY20121}. The proof of the
desired bounds of the entries of $\Sigma^{1/2}\mathcal{G}\Sigma^{1/2}$
will be postponed to the next section, which can be viewed as our main
technical ingredient of this paper.

\begin{thm}[(Green function comparison theorem around $\lambda_r$)]
\label{Green} Let $\mathcal{W}^{\mathbf{u}}$ and $\mathcal{W}^{\mathbf{v}}$
be two sample covariance matrices in Theorem~\ref{thm.0.0}. Let $F$ be
a real function satisfying
%
\begin{equation}\label{1026.20}
\sup_x\bigl|F^{(k)}(x)\bigr|/\bigl(|x|+1\bigr)^{C}\leq C,\qquad
k=0,1,2,3,4
\end{equation}
for some positive constant $C$. There exist $\varepsilon_0>0$ and
$N_0\in\mathbb{N}$, such that for any positive constant $\varepsilon
<\varepsilon_0$, $N\geq N_0$ and for any\vspace*{1.5pt} real numbers $E,E_1$ and $E_2$
satisfying $|E-\lambda_r|,|E_1-\lambda_r|,|E_2-\lambda_r|\leq
N^{-2/3+\varepsilon}$
and $\eta=N^{-2/3-\varepsilon}$, we have for $z=E+\mathbf{i}\eta$,
%
\begin{equation}\label{11.11.00}
\bigl\llvert \mathbb{E}F\bigl(N\eta\Im m^{\mathbf{u}}_N(z)
\bigr)-\mathbb{E}F\bigl(N\eta\Im m^{\mathbf{v}}_N(z)\bigr)\bigr
\rrvert \leq N^{-C^\prime\varepsilon}
\end{equation}
and
%
\begin{eqnarray}\label{11.11.01}
&&\biggl\llvert \mathbb{E}F \biggl(N\int_{E_1}^{E_2}
\Im m_N^{\mathbf{u}}(x+\mathbf {i}\eta)\,dx \biggr)-\mathbb{E}F
\biggl(N\int_{E_1}^{E_2} \Im m_N^{\mathbf{v}}(x+
\mathbf{i}\eta)\,dx \biggr)\biggr\rrvert 
\nonumber
\\[-8pt]
\\[-8pt]
\nonumber
&&\qquad\leq N^{-C^\prime
\varepsilon}
\end{eqnarray}
with some positive constant $C^\prime$ if either $\mathbf{A}$ or
$\mathbf{B}$ in Theorem~\ref{thm.0.0} holds.
\end{thm}

Now we are at the stage to prove our main results assuming Theorem~\ref{Green}.
\begin{pf*}{Proof of Theorem~\ref{thm.0.0}}
Given Theorems \ref{lem.2.2}--\ref{thm.1206.7} and \ref{Green}, the proof of Theorem~\ref{thm.0.0} is nearly the same as that
for the null case in \cite{PY2012}
(see the proof of Theorem~1.1 therein). Due to the similarity, here we
only sketch the main route and leave the details to the reader. We
start from Theorem~\ref{thm.1206.7}(i), which states that for any
$\zeta>0$, there exists some positive constant $C_\zeta$ such that
$|\lambda_1(\mathcal{W})-\lambda_r|\leq N^{-2/3}\varphi^{C_\zeta}$ with
$\zeta$-high probability. Hence, it suffices to verify\vspace*{1pt} (\ref{0.1.2.3.2.1}) for $s\in[-\frac{3}2\varphi^{C_\zeta}, \frac{3}2\varphi
^{C_\zeta}]$. To this end, we denote $E_\zeta=\lambda_r+2N^{-2/3}\varphi
^{C_\zeta}$ and set $E=\lambda_r+sN^{-2/3}$.
With the above restriction on $s$, one can always assume that $E\leq
E_\zeta-\frac{1}2 N^{-2/3}\varphi^{C_\zeta}$. Denoting $\eta
_1=N^{-2/3-9\varepsilon_1}$ and $\ell=\frac{1}2 N^{-2/3-\varepsilon_1}$
with any given small constant $\varepsilon_1>0$, we record the
following inequality from Corollary~5.1 of \cite{BPZ2013}:
%
\begin{eqnarray}
&& \mathbb{E}h\biggl(\frac{N}{\pi}\int_{E-\ell}^{E_\zeta}
\Im m_N(y+\mathbf{i}\eta _1)\,dy\biggr)\nonumber\\
&& \label{2014050501}\qquad\leq  \mathbb{P}
\bigl(\lambda_1(\mathcal{W})\leq E\bigr)
\\
\nonumber
&&\qquad\leq \mathbb{E}h\biggl(\frac{N}{\pi}\int_{E+\ell}^{E_\zeta}
\Im m_N(y+\mathbf {i}\eta_1)\,dy\biggr)+O\bigl(\exp
\bigl(-\varphi^{C_\zeta}\bigr)\bigr),
\end{eqnarray}
where $h$ is a smooth cutoff function satisfying the condition of $F$
in Theorem~\ref{Green}; see Corollary~5.1 of \cite{BPZ2013} for the
definition of the function $h$. (\ref{2014050501}) states that $\mathbb{P}(\lambda_1(\mathcal{W})\leq E)$ can be squeezed by the expectations
of two functionals of the Stieltjes transform.
Now, setting $\varepsilon=9\varepsilon_1$, $\eta=\eta_1$, $F(x)=h(x/\pi
)$, $E_1=E-\ell$ and $E_2=E_\zeta$ in (\ref{11.11.01}) we obtain
%
\begin{eqnarray}
&& \mathbb{E}h \biggl(\frac{N}{\pi}\int_{E-\ell}^{E_\zeta}
\Im m_N^{\mathbf
{u}}(x+\mathbf{i}\eta_1)\,dx \biggr)
\nonumber
\\[-8pt]
\label{2014050502}\\[-8pt]
\nonumber
&& \qquad \leq\mathbb{E}h \biggl(\frac{N}{\pi}\int_{E-\ell}^{E_\zeta}
\Im m_N^{\mathbf{v}}(x+\mathbf{i}\eta_1)\,dx \biggr)+
\frac{1}2 N^{-\delta},
\end{eqnarray}
for sufficiently large $N$, where we took $\delta=\frac{1}2 C^\prime
\varepsilon$ (say).
Employing the second inequality in (\ref{2014050501}) via replacing $E$
by $E-2\ell$, we also have
\[
\mathbb{P}\bigl(\lambda_1\bigl(\mathcal{W}^{\mathbf{u}}\bigr)
\leq E-2\ell\bigr)\leq\mathbb{E}h\biggl(\frac{N}{\pi}\int
_{E-\ell}^{E_\zeta}\Im m_N^{\mathbf{u}}(y+
\mathbf {i}\eta_1)\,dy\biggr)+O\bigl(\exp\bigl(-\varphi^{C_\zeta}
\bigr)\bigr),
\]
which together with (\ref{2014050502}) implies that for sufficiently
large $N$,
\[
\mathbb{P}\bigl(\lambda_1\bigl(\mathcal{W}^{\mathbf{u}}\bigr)
\leq E-2\ell\bigr)\leq\mathbb {E}F \biggl(\frac{N}{\pi}\int
_{E-\ell}^{E_\zeta} \Im m_N^{\mathbf{v}}(x+
\mathbf{i}\eta_1)\,dx \biggr)+N^{-\delta}.
\]
Using the first inequality in (\ref{2014050501}) yields
%
\begin{equation}\label{20140505031}
\mathbb{P}\bigl(\lambda_1\bigl(\mathcal{W}^{\mathbf{u}}\bigr)
\leq E-2\ell\bigr)\leq\mathbb {P}\bigl(\lambda_1\bigl(
\mathcal{W}^{\mathbf{v}}\bigr)\leq E\bigr)+N^{-\delta}.
\end{equation}
Switching the roles of $\mathbf{u}$ and $\mathbf{v}$, we can
analogously derive that
%
\begin{equation}\label{20140505041}
\mathbb{P}\bigl(\lambda_1\bigl(\mathcal{W}^{\mathbf{v}}\bigr)
\leq E\bigr)\leq\mathbb {P}\bigl(\lambda_1\bigl(
\mathcal{W}^{\mathbf{u}}\bigr)\leq E+2\ell\bigr)+N^{-\delta}.
\end{equation}
(\ref{20140505031}) and (\ref{20140505041}) then lead to (\ref
{0.1.2.3.2.1}). Hence, we conclude the sketch of the proof.
\end{pf*}
\begin{pf*}{Proof of Corollary~\ref{thm.1.1}}
Corollary~\ref{thm.1.1}
follows from Theorem~\ref{thm.0.0}, Theorem~1 of \cite{Karoui2007} and
Proposition~2 of \cite{Onatski2008} immediately.
\end{pf*}
Now, before commencing the proof of Corollary~\ref{thm.1205.1}, we
record the following lemma whose proof will be provided in the
supplementary material \cite{BPZ2014}.

\begin{lem} \label{lem.2014050501}
Assume that $\Sigma$ satisfies the
assumption of Corollary~\ref{thm.1205.1}. Then $\Sigma$ also satisfies
Condition \ref{con.1.1}\textup{(iii)}. In addition, we have (\ref{1206.12}).
\end{lem}

\begin{pf*}{Proof of Corollary~\ref{thm.1205.1}}
With the aid of
Lemma~\ref{lem.2014050501}, we see that $\mathcal{W}$ satisfies
Condition \ref{con.1.1} thus Theorem~\ref{thm.0.0} can be adopted. Now
we invoke the fact that the real Wishart matrix with population
covariance matrix $\operatorname{diag}(\lambda_{1}(\Sigma),\ldots, \lambda
_{M}(\Sigma))$ satisfy the conditions of Theorem~1.6 of \cite{FP2009}.
Moreover, taking the property of orthogonal invariance for Gaussian
matrices into account, we know the result of \cite{FP2009} holds for
all Wishart matrices with population covariance matrix $\Sigma$
(possibly not diagonal) whose eigenvalues satisfy the condition in
Corollary~\ref{thm.1205.1}. We remind here the parameters $N$ and $p$
in \cite{FP2009} are corresponding to our $M$ and $N$, respectively.
Hence, with (\ref{1206.12}) at hand, choosing the Wishart matrix with
population covariance matrix $\Sigma$ as the reference matrix and
combining our Theorem~\ref{thm.0.0} with Theorem~1.6 of \cite{FP2009},
we can complete the proof.
\end{pf*}
It remains to prove Theorem~\ref{Green} in this section.
\begin{pf*}{Proof of Theorem~\ref{Green}}
To simplify the presentation, we will only verify (\ref{11.11.00}) in
detail below. The proof of (\ref{11.11.01}) can be taken similarly,
thus we just leave it to the reader. As a compensation, some necessary
modifications for the proof of (\ref{11.11.01}) will be highlighted in
Remarks \ref{rem.1207.100} and \ref{rem.20140501}.
Now, let $\gamma\in\{1,2,\ldots,N+1\}$ and set $X_\gamma$ to be the
matrix whose first $\gamma-1$ columns are the same as those of
$X^{\mathbf{v}}$ and the remaining $N-\gamma+1$ columns are the same as
those of $X^{\mathbf{u}}$. Especially, we have $X_1=X^{\mathbf{u}}$ and
$X_{N+1}=X^{\mathbf{v}}$. Correspondingly, we set
\[
W_{N,\gamma}=X_\gamma^{*}\Sigma X_{\gamma},\qquad
\mathcal{W}_{N,\gamma
}=\Sigma^{1/2}X_{\gamma}X^{*}_{\gamma}
\Sigma^{1/2}.
\]
Then (\ref{11.11.00}) can be achieved via checking that for every
$\gamma$ the following estimate holds:
\[
\mathbb{E} F \bigl(\eta\Im \Trr (W_{N,{\gamma}}-z)^{-1} \bigr)-
\mathbb{E} F \bigl(\eta\Im \Trr (W_{N,{\gamma+1}}-z)^{-1} \bigr)=O
\bigl(N^{-1-C'\varepsilon}\bigr).
\]
Observing the fact that $X_{\gamma}$ and $X_{\gamma+1}$ only differ in
the $\gamma$th column yields $X_{\gamma}^{(\gamma)}=X_{\gamma
+1}^{(\gamma)}$,
which directly implies $W_{N,\gamma}^{(\gamma)}=W_{N,\gamma+1}^{(\gamma
)}$ and $\mathcal{W}_{N,\gamma}^{(\gamma)}=\mathcal{W}_{N,\gamma
+1}^{(\gamma)}$. Therefore, we can also write
%
\begin{eqnarray}
&&\mathbb{E} F \bigl(\eta\Im \Trr (W_{N,{\gamma}}-z)^{-1} \bigr)-
\mathbb{E} F \bigl(\eta\Im \Trr (W_{N,{\gamma+1}}-z)^{-1} \bigr)
\nonumber
\\
&&\qquad= \bigl(\mathbb{E} F \bigl(\eta\Im \Trr (W_{N,{\gamma}}-z)^{-1}
\bigr)-\mathbb{E} F \bigl(\eta\Im\bigl[\Trr \bigl(W_{N,{\gamma}}^{(\gamma
)}-z
\bigr)^{-1}-z^{-1}\bigr] \bigr) \bigr)
\nonumber
\\[-8pt]
\label{k.20}
\\[-8pt]
\nonumber
&&\qquad\quad{}- \bigl(\mathbb{E} F \bigl(\eta\Im \Trr (W_{N,{\gamma+1}}-z)^{-1}
\bigr)\\
&& \hspace*{17pt}\qquad\quad{}-\mathbb{E} F \bigl(\eta\Im\bigl[\Trr \bigl(W_{N,{\gamma+1}}^{(\gamma
)}-z
\bigr)^{-1}-z^{-1}\bigr] \bigr) \bigr).\nonumber
\end{eqnarray}
Since the comparison process will greatly rely on the moment matching
condition, it will be more convenient to work with the following set:
\[
\mathcal{M}_{k}(i):=\bigl\{\bigl(\mathbb{E}(\Re\sqrt{N}x_{ji})^l(
\Im\sqrt {N}x_{ji})^m,j,l,m\bigr)\dvtx j=1,\ldots, M, m+l\leq k
\bigr\},
\]
that is, the set of all moments up to order $k$ of the entries of $\sqrt
{N}\mathbf{x}_i$, where its elements are indexed by $j,l,m$. In the
spirit of (\ref{k.20}), it suffices to show, for any sample covariance
matrix $W_N$ satisfying Condition \ref{con.1.1}, the following Lemmas
\ref{lem.11.11.00} and \ref{lem.y.4.6} hold.

\begin{lem} \label{lem.11.11.00} Let $F$ be a real function satisfying
(\ref{1026.20}) and $z=E+\mathbf{i}\eta$. For any random matrix $W_N$
satisfying Condition \ref{con.1.1}, if
$|E-\lambda_r|\leq N^{-2/3+\varepsilon}$ and $N^{-2/3-\varepsilon}\leq
\eta\ll N^{-2/3}$
for some $\varepsilon>0$, there exists some positive constant $C$
independent of $\varepsilon$ such that
%
\begin{eqnarray}
&& \mathbb{E} F\bigl(N\eta\Im m_N(z)\bigr)-\mathbb{E}F\bigl(N\eta\Im
\bigl[m_N^{(i)}(z)-(Nz)^{-1}\bigr]\bigr)
\nonumber
\\[-8pt]
\label{11.11.100}\\[-8pt]
\nonumber
&&\qquad =A
\bigl(X^{(i)}, \mathcal {M}_2(i)\bigr)+N^{-1-C\varepsilon}
\end{eqnarray}
when $\Sigma$ is diagonal, and
%
\begin{eqnarray}
&& \mathbb{E} F\bigl(N\eta\Im m_N(z)\bigr)-\mathbb{E}F\bigl(N\eta\Im
\bigl[m_N^{(i)}(z)-(Nz)^{-1}\bigr]\bigr)
\nonumber
\\[-8pt]
\label{11.11.101}\\[-8pt]
\nonumber
&&\qquad =B
\bigl(X^{(i)}, \mathcal {M}_4(i)\bigr)+N^{-1-C\varepsilon}
\end{eqnarray}
for general $\Sigma$,
where $A(X^{(i)},\mathcal{M}_2(i))$ is a functional depending on the
distribution of $X^{(i)}$ and $\mathcal{M}_2(i)$ only and similarly
$B(X^{(i)},\mathcal{M}_4(i))$ is
a functional depending on the distribution of $X^{(i)}$ and $\mathcal
{M}_4(i)$ only.
\end{lem}

\begin{rem} \label{rem.1207.100}
To verify (\ref{11.11.01}), actually we need to show two equalities
analogous to (\ref{11.11.100}) and (\ref{11.11.101}), obtained via
replacing
\[
F\bigl(N\eta\Im m_N(z) \mbox{ and }F\bigl(N\eta\Im
\bigl[m_N^{(i)}(z)-(Nz)^{-1}\bigr]\bigr)\bigr)
\]
by
\[
F\biggl(N\int_{E_1}^{E_2}\Im m_N(x+\mathbf{i}\eta)\,dx\biggr)\quad \mbox{and}\quad
F\biggl(N\int_{E_1}^{E_2}\Im\bigl[m_N^{(i)}(x+\mathbf{i}\eta)-(Nz)^{-1}\bigr]\,dx\biggr),
\]
respectively, in (\ref{11.11.100}) and (\ref{11.11.101}) and
correspondingly replacing $A(\cdot,\cdot)$ and $B(\cdot,\cdot)$ by some
other functionals $\tilde A(\cdot,\cdot)$ and $\tilde B(\cdot,\cdot)$.
\end{rem}

Now, to differentiate, we denote the set $\mathcal{M}_k(i)$ for
$X^{\mathbf{u}}$ and $X^{\mathbf{v}}$ by $\mathcal{M}_k^{\mathbf
{u}}(i)$ and $\mathcal{M}_k^{\mathbf{v}}(i)$, respectively. Then, we
also have the following lemma.

\begin{lem} \label{lem.y.4.6}
Under Condition \ref{con.1.1} and the
assumptions in Lemma~\ref{lem.11.11.00}, there exist some positive
constants $c_0$ and $C>1$, such that the following statements hold. If
$\mathcal{W}^{\mathbf{u}}$ matches $\mathcal{W}^{\mathbf{v}}$ to order
$2$, we have
%
\begin{equation}\label{p.p.p.1.1.1}
\max_{\gamma}\bigl|A\bigl(X_\gamma^{(\gamma)},
\mathcal{M}^{\mathbf{u}}_2(\gamma )\bigr)-A\bigl(X_\gamma^{(\gamma)},
\mathcal{M}^{\mathbf{v}}_2(\gamma )\bigr)\bigr|=O\bigl(e^{-c_0(\log N)^C}
\bigr).
\end{equation}
If $\mathcal{W}^{\mathbf{u}}$ matches $\mathcal{W}^{\mathbf{v}}$ to
order $4$, we have
%
\begin{equation}\label{p.p.p.2.2.2}
\max_{\gamma}\bigl|B\bigl(X_\gamma^{(\gamma)},
\mathcal{M}^{\mathbf{u}}_4(\gamma )\bigr)-B\bigl(X_\gamma^{(\gamma)},
\mathcal{M}^{\mathbf{v}}_4(\gamma )\bigr)\bigr|=O\bigl(e^{-c_0(\log N)^C}
\bigr).
\end{equation}
Here, $A(\cdot,\cdot)$ and $B(\cdot,\cdot)$ are the functionals in
Lemma~\ref{lem.11.11.00}.
\end{lem}

It is obvious that (\ref{11.11.00}) follows from Lemma~\ref
{lem.11.11.00} and Lemma~\ref{lem.y.4.6}. The proof of (\ref{11.11.01})
can be taken analogously. Thus, we conclude the proof of Theorem~\ref
{Green} assuming the validity of Lemmas \ref{lem.11.11.00} and \ref{lem.y.4.6}.
\end{pf*}
We leave the proof of Lemma~\ref{lem.y.4.6} to the supplementary
material \cite{BPZ2014} and prove Lemma~\ref{lem.11.11.00} in the
sequel. Without loss of generality, we will just check the statements
in Lemma~\ref{lem.11.11.00} for $i=1$. The others are just analogous.
To commence the proof, we will need the following lemma as an input,
whose proof will also appear in the supplementary material \cite{BPZ2014}.

\begin{lem} \label{lem.20140505123}
Under the assumptions on $z$ and $F$
in Lemma~\ref{lem.11.11.00}, for any given $\zeta>0$, there exists some
positive constant $C$, such that
%
\begin{eqnarray}
&&\hspace*{4pt} F\bigl(N\eta\Im m_N(z)\bigr)-F\bigl(N\eta\Im
\bigl[m_N^{(1)}(z)-(Nz)^{-1}\bigr]\bigr)
\nonumber
\\[-8pt]
\label{515.1}\\[-8pt]
\nonumber
&&\qquad=\sum_{k=1}^{3}\frac{1}{k!}F^{(k)}
\bigl(N\eta\Im \bigl[m_N^{(1)}(z)-(Nz)^{-1}
\bigr]\bigr) (\Im y)^k+O\bigl(N^{-4/3+C\varepsilon}\bigr)
\end{eqnarray}
holds with $\zeta$-high probability, where
%
\begin{equation}\label{11.11.30}
y:=\eta zG_{11}\mathbf{r}_1^*\bigl(
\mathcal{G}^{(1)}\bigr)^2\mathbf{r}_1.
\end{equation}
Moreover, we have
%
\begin{equation}\label{515.299999}
\bigl|\mathbf{r}_1^*\bigl(\mathcal{G}^{(1)}
\bigr)^2\mathbf{r}_1\bigr|\leq N^{1/3+C\varepsilon},\qquad |y|\leq
N^{-1/3+C\varepsilon}
\end{equation}
with $\zeta$-high probability.
\end{lem}

With Lemma~\ref{lem.20140505123}, we now start to prove Lemma~\ref
{lem.11.11.00} for $i=1$.

\begin{pf*}{Proof of Lemma~\ref{lem.11.11.00} (for $i=1$)}
Now, starting from (\ref{11.11.30}), we further decompose $y$ and then
pick out the leading terms in the decomposition. Specifically, we set
%
\begin{equation}\label{1207.1}
D:=\frac{m_0-G_{11}}{G_{11}}=-m_0\cdot\bigl(z+z\mathbf{r}_1^*
\mathcal {G}^{(1)}\mathbf{r}_1\bigr)-1,
\end{equation}
which is implied by the Schur complement $G_{11}=-1/(z+z\mathbf
{r}_1^*\mathcal{G}^{(1)}\mathbf{r}_1)$. At first, by (i) of Theorem~\ref{lem.2.2} and (ii) of Theorem~\ref{thm.3.90} we can see that
$G_{ii}(z)\sim1$ with $\zeta$-high probability. Moreover, with $\zeta$-high probability we can write
%
\begin{equation}\label{11.11.20}
G_{11}=\frac{m_0}{D+1}=m_0\sum
_{k=0}^{\infty}(-D)^{k}
\end{equation}
since $|D|\leq N^{-1/3+C\varepsilon}$ for some positive constant $C$,
which is implied by Theorem~\ref{lem.2.2}(ii) and Theorem~\ref{thm.3.90}(ii).
Inserting (\ref{11.11.20}) into (\ref{11.11.30}), we can write
%
\begin{equation}\label{201405130101}
y=\sum_{k=1}^{\infty}y_k,\qquad
y_k:=\eta z m_0(-D)^{k-1}\mathbf
{r}_1^{*}\bigl(\mathcal{G}^{(1)}
\bigr)^2\mathbf{r}_1.
\end{equation}
By (\ref{515.299999}) and the aforementioned bound for $D$, we can
easily get
%
\begin{equation}\label{201404130102}
|y_k|=O\bigl(N^{-k/3+C\varepsilon}\bigr)
\end{equation}
with $\zeta$-high probability, which directly implies that
%
\begin{eqnarray}
\Im y &=& \Im y_1+\Im y_2+\Im y_3+O
\bigl(N^{-4/3+C\varepsilon}\bigr),
\nonumber
\\[-2pt]
\label{22220000111010101010}
(\Im y)^2 &=& (\Im y_1)^2+2\Im
y_1\Im y_2+ O\bigl(N^{-4/3+C\varepsilon}\bigr),\\[-2pt]
\nonumber
(\Im
y)^3 &=& (\Im y_1)^3+O\bigl(N^{-4/3+C\varepsilon}
\bigr)
\end{eqnarray}
hold with $\zeta$-high probability. By the discussions in the proof of
Lemma~\ref{lem.20140505123} in the supplementary material\vspace*{1pt} \cite
{BPZ2014}, one can see that $N\eta \Im m_N(z)=O(N^{C\varepsilon})$ and
$N\eta\Im[m_N^{(1)}(z)-(Nz)^{-1}]=O(N^{C\varepsilon})$ with $\zeta$-high probability for any given $\zeta>0$. Consequently, in light of
the assumption\vspace*{1pt} (\ref{1026.20}), we see that for any real number $t_N$
between $N\eta \Im m_N(z)$ and $N\eta\Im[m_N^{(1)}(z)-(Nz)^{-1}]$,
there exist
%
\begin{equation}\label{2014051301010101}
\bigl|F^{(k)}(t_N)\bigr|=O\bigl(N^{C\varepsilon}\bigr), \qquad k=0,1,2,3,4
\end{equation}
with $\zeta$-high probability. Moreover, we have the deterministic
bound $|m_N(z)|,\break  |m_N^{(1)}(z)|=O(\eta^{-1})$, which implies $|N\eta\Im
m_N(z)|,|N\eta\Im[m_N^{(1)}(z)-(Nz)^{-1}]|=O(N)$. Thus, using (\ref
{1026.20}) again we have the deterministic bound
$|F^{(k)}(t_N)|=O(N^C), k=0,1,2,3,4$, for any real number $t_N$ between
$N\eta \Im m_N(z)$ and $N\eta\Im[m_N^{(1)}(z)-(Nz)^{-1}]$.
Analogously, by using the fact that $\Vert\mathcal{G}^{-1}\Vert_{\mathrm{op}}=O(\eta
^{-1})$ and condition (\ref{1206.5}), we can get the deterministic
bound $|y|=O(N^{C})$, $|y_k|=O(N^{C(k)})$ with some positive constants
$C$ and $C(k)$ (depending on $k$), for $k=0,1,2,3,4$. Then by (\ref
{515.1}), (\ref{201405130101})--(\ref{2014051301010101}) and the
deterministic bounds above, it is not difficult to find that
%
\begin{eqnarray}
&&\mathbb{E}F\bigl(N\eta\Im m_N(z)\bigr)-\mathbb{E}F\bigl(N\eta\Im
\bigl[m_N^{(1)}(z)-(Nz)^{-1}\bigr]\bigr)
\nonumber
\\[-2pt]
&&\qquad = \mathbb{E}F^{(1)}\bigl(N\eta\Im\bigl[m_N^{(1)}(z)-(Nz)^{-1}
\bigr]\bigr) (\Im y_1+\Im y_2+\Im y_3)
\nonumber
\\[-10pt]
\label{11.12.13}
\\[-10pt]
\nonumber
&&\qquad \quad{}+\mathbb{E}F^{(2)}\bigl(N\eta\Im\bigl[m_N^{(1)}(z)-(Nz)^{-1}
\bigr]\bigr) \bigl(\tfrac{1}2(\Im y_1)^2+\Im
y_1\Im y_2\bigr)
\nonumber
\\[-2pt]
&&\qquad\quad{}+\mathbb{E}F^{(3)}\bigl(N\eta\Im\bigl[m_N^{(1)}(z)-(Nz)^{-1}
\bigr]\bigr) \bigl(\tfrac{1}6\Im y_1\bigr)^3+O
\bigl(N^{-4/3+C\varepsilon}\bigr).\nonumber
\end{eqnarray}
Toward the right-hand side of (\ref{11.12.13}), our task is to extract
the terms depending on $X^{(1)}$ and $\mathcal{M}_k(1)$ ($k=2$ or $4$)
only and bound the remaining terms. For the latter, we will need the
following crucial lemma on bounding the\vspace*{-2pt} entries of $\Sigma^{1/2}\mathcal
{G}^{(1)}\Sigma^{1/2}$.

\begin{lem} \label{lem.11.11.105} Let $z=E+\mathbf{i}\eta$ with $
|E-\lambda_r|\leq N^{-2/3+\varepsilon}$ and $N^{-2/3-\varepsilon}\leq
\eta\ll N^{-2/3}$
for some $\varepsilon>0$. When $\Sigma$ is diagonal, for any given
$\zeta>0$, we have
%
\begin{eqnarray}
 \bigl|\bigl(\mathcal{G}^{(1)}(z)\bigr)_{ij}\bigr| &\leq &
N^{C\varepsilon}\quad \mbox{and}
\nonumber
\\[-10pt]
\label{502.20}
\\[-10pt]
\nonumber
\bigl|\bigl(\bigl[\mathcal{G}^{(1)}(z)
\bigr]^2\bigr)_{ij}\bigr| &\leq &  N^{1/3+C\varepsilon},\qquad i,j\in \{1,
\ldots,M\}
\end{eqnarray}
hold with $\zeta$-high probability for some positive constant $C$
independent of $\varepsilon$.

For general $\Sigma$ satisfying Condition \ref{con.1.1}\textup{(iii)}, we have
for any given $\zeta>0$,
%
\begin{equation}\label{502.21}
\bigl|\bigl(\Sigma^{1/2}\mathcal{G}^{(1)}(z)\Sigma^{1/2}
\bigr)_{kk}\bigr|\leq N^{1/3+C\varepsilon},\qquad k\in\{1,\ldots,M\}
\end{equation}
hold with $\zeta$-high probability for some positive constant $C$
independent of $\varepsilon$.
\end{lem}

\begin{rem}\label{rem.20140501}
When we prove (\ref{11.11.01}), as
mentioned above, we actually need to verify the statement in Remark~\ref
{rem.1207.100}. To this end, we need to strengthen the bounds in (\ref
{502.20}) and (\ref{502.21}) to hold with $\zeta$-high probability
uniformly on the set $\{z=E+\mathbf{i}\eta\dvtx  |E-\lambda_r|\leq
N^{-2/3+\varepsilon}\mbox{ and } N^{-2/3-\varepsilon}\leq\eta\ll
N^{-2/3}\}$.
These uniform bounds are necessary for the proof of the statement in
Remark~\ref{rem.1207.100}, since some integrations taken w.r.t. the
real part of $z$ are involved in the discussion. These stronger bounds
can be obtained from the bounds for single point in (\ref{502.20}) and (\ref{502.21}) through some routine $\varepsilon$-net and
Lipschitz continuity argument. One can refer to the extension from (\ref
{x.x.x}) to (\ref{1207.30}) below for a similar argument.
\end{rem}

Lemma~\ref{lem.11.11.105} is our main technical task whose proof will
be postponed to the next section separately. Now, with Lemma~\ref
{lem.11.11.105} granted, we prove (\ref{11.11.100}) and (\ref
{11.11.101}) in the sequel.
At first, we will verify (\ref{11.11.100}) for diagonal $\Sigma$. We
start with the third term on the r.h.s. of (\ref{11.12.13}). Denoting
$\varpi:=\Im\eta zm_0$ and $\varrho:=\Re\eta zm_0$, we have
\[
\Im y_1=\varpi\bigl(\Re\mathbf{r}_1^*\bigl(
\mathcal{G}^{(1)}\bigr)^2\mathbf {r}_1\bigr)+
\varrho\bigl(\Im\mathbf{r}_1^*\bigl(\mathcal{G}^{(1)}
\bigr)^2\mathbf{r}_1\bigr).
\]
To further simplify the exposition, we denote the real part and
imaginary part of a complex number $A$ by $A[0]$ and $A[1]$
respectively. Introducing\vspace*{1pt} the notation $\mathbb{E}_i$ to denote the
expectation with respect to $\mathbf{x}_i$, we see that
$\mathbb{E}_1(\Im y_1)^3$ is a summation of finite terms of the form
%
\begin{equation}\label{12.12.12}
\hspace*{15pt}\mathbf{1}_{\{a+b=3\}}\varpi^{a}\varrho^{b} \sum
_{k_1,\ldots,k_6}\prod_{i=1}^3
\bigl(\Sigma^{1/2}\bigl(\mathcal{G}^{(1)}\bigr)^2
\Sigma^{1/2}\bigr)_{k_{2i-1}k_{2i}}[\alpha _i]\mathbb{E}
\prod_{l=1}^6x_{k_l,1}[
\beta_l],
\end{equation}
where $\alpha_i,\beta_l$ are $0$ or $1$ and $a,b$ are nonnegative
integers. Hence, it suffices to analyze the quantities of the form (\ref
{12.12.12}) below.

We classify the terms in the summation (\ref{12.12.12}) by various
coincidence conditions of the indices $k_1,\ldots,k_6$. If there is a
$k_j$ appearing only once in $\{k_1,\ldots,k_6\}$, then this term is
zero obviously, due to the independence and centering of the entries of
$X$. Now we proceed to those terms in which each $k_j$ appears exactly
twice. Apparently, these terms only depend on $X^{(1)}$ and $\mathcal
{M}_2(1)$. Finally, it remains to consider the terms that there is at
least one $k_j$
appearing at least three times and no $k_j$ appearing\vspace*{1pt} only once. It is
obviously that the total number of such terms is $O(N^2)$. Putting this
observation and (\ref{502.20}) in Lemma~\ref{lem.11.11.105} together
yields the fact that
the total contribution of these terms is less than
\begin{eqnarray*}
 CN^{-1}|\eta zm_0|^3\max
_{ij}\bigl|\bigl(\Sigma^{1/2}\bigl(\mathcal{G}^{(1)}
\bigr)^2\Sigma ^{1/2}\bigr)_{ij}\bigr|^3
&\leq &  CN^{-1}|\eta zm_0|^3\max
_{ij}\bigl|\bigl(\bigl(\mathcal {G}^{(1)}\bigr)^2\bigr)_{ij}\bigr|^3\\
&=& O
\bigl(N^{-2+C\varepsilon}\bigr)
\end{eqnarray*}
with $\zeta$-high probability. Since $|(\Sigma^{1/2}\mathcal
{G}^{(1)}\Sigma^{1/2})_{ij}|$ are trivially bounded by $O(\eta^{-1})$,
one can see that the above bound also holds in expectation by the
definition of $\zeta$-high probability.
Therefore, we deduce that
%
\begin{equation}\label{12.12.100}
\mathbb{E}(\Im y_1)^3=A_1
\bigl(X^{(1)},\mathcal{M}_2(1)\bigr)+O\bigl(N^{-2+C\varepsilon
}
\bigr)
\end{equation}
for some functional $A_1$ depending on the distribution of $X^{(1)}$
and $\mathcal{M}_2(1)$ only.

Now, for the first and second term on the right-hand side of (\ref
{11.12.13}) we can deal with them analogously. Note that by (\ref
{1207.1}) and the definitions of $y_2,y_3$, one can see that
\begin{eqnarray*}
y_2 &=& \eta z^2m_0^2
\mathbf{r}_1^*\mathcal{G}^{(1)}\mathbf{r}_1
\cdot \mathbf{r}_1^*\bigl(\mathcal{G}^{(1)}
\bigr)^2\mathbf{r}_1 +\eta zm_0
\cdot(1+zm_0)\cdot\mathbf{r}_1^*\bigl(
\mathcal{G}^{(1)}\bigr)^2\mathbf {r}_1,
\\
y_3 &=& \eta z^3m_0^3\bigl(
\mathbf{r}_1^*\mathcal{G}^{(1)}\mathbf{r}_1
\bigr)^2\cdot \mathbf{r}_1^*\bigl(\mathcal{G}^{(1)}
\bigr)^2\mathbf{r}_1\\
&&{} +2\eta z^2m_0
\cdot(1+zm_0)\cdot\mathbf{r}_1^*\mathcal{G}^{(1)}
\mathbf {r}_1\cdot\mathbf{r}_1^*\bigl(
\mathcal{G}^{(1)}\bigr)^2\mathbf{r}_1
\\
&&{}+\eta z m_0\cdot(1+zm_0)\mathbf{r}^*_1
\mathcal {G}^{(1)}\mathbf{r}_1.
\end{eqnarray*}
Expanding each term above, then by a routine but detailed discussion on
the coincidence condition of the indices as what we have done to the
third term on the right-hand side of (\ref{11.12.13}) above, we can
actually get
%
\begin{equation}\label{12.12.101}
\mathbb{E} \bigl((\Im y_1)^2+2(\Im y_1)
(\Im y_2) \bigr)=A_2\bigl(X^{(1)},
\mathcal{M}_2(1)\bigr)+O\bigl(N^{-5/3+C\varepsilon}\bigr)
\end{equation}
and
%
\begin{equation}\label{12.12.102}
\mathbb{E}(\Im y_1+\Im y_2+\Im y_3)=A_3
\bigl(X^{(1)},\mathcal {M}_2(1)\bigr)+O\bigl(N^{-4/3+C\varepsilon}
\bigr)
\end{equation}
for some functionals $A_2$ and $A_3$ depending on $X^{(1)}$ and
$\mathcal{M}_2(1)$ only.
Inserting (\ref{12.12.100})--(\ref{12.12.102}) into (\ref{11.12.13}),
we obtain (\ref{11.11.100}).

Now, we go ahead to investigate (\ref{11.11.101}) for more general
$\Sigma$. At first, we revisit the canonical form of the terms in the
expansion of $\mathbb{E}_1(\Im y_1)^3$, that is, (\ref{12.12.12}). Note
that for (\ref{12.12.12}), it suffices to bound the terms in which all
$k_i,i=1,\ldots,6$ are the same, since all the other terms only depend
on the distribution of $X^{(1)}$ and $\mathcal{M}_4(1)$. In other
words, we need to bound the terms in which $\mathbb{E}|x_{k1}|^6$
appears. Analogously, writing $\mathbb{E}_1 ((\Im y_1)^2+2(\Im
y_1)(\Im y_2) )$ and $\mathbb{E}_1(\Im y_1+\Im y_2+\Im y_3)$ as
some summations of terms in the forms similar to (\ref{12.12.12}),
again, we only need to address the terms containing $\mathbb
{E}|x_{k1}|^6$ as a factor. It is not difficult to see after simple
calculations that the total contribution of such terms in $\mathbb
{E}_1(\Im y_1)^3$, $\mathbb{E}_1 ((\Im y_1)^2+2(\Im y_1)(\Im
y_2) )$ and $\mathbb{E}_1(\Im y_1+\Im y_2+\Im y_3)$ can be bounded by
%
\begin{eqnarray}
&& CN^{-3}\sum_k\eta^3\bigl|\bigl(
\Sigma^{1/2}\bigl(\mathcal{G}^{(1)}\bigr)^2
\Sigma^{1/2}\bigr)_{kk}\bigr|^3 \nonumber\\
&&\label{502.45}\qquad{}+ CN^{-3}
\sum_{k}\eta^2\bigl|\bigl(
\Sigma^{1/2}\bigl(\mathcal{G}^{(1)}\bigr)^2\Sigma
^{1/2}\bigr)_{kk}\bigr|\bigl|\bigl(\Sigma^{1/2}\bigl(
\mathcal{G}^{(1)}\bigr)^2\Sigma^{1/2}
\bigr)_{kk}\bigr|^2
\\
\nonumber
&&\qquad{}+ CN^{-3}\sum_{k}\eta\bigl|\bigl(
\Sigma^{1/2}\bigl(\mathcal{G}^{(1)}\bigr)\Sigma ^{1/2}
\bigr)_{kk}\bigr|^2\bigl|\bigl(\Sigma^{1/2}\bigl(
\mathcal{G}^{(1)}\bigr)^2\Sigma^{1/2}
\bigr)_{kk}\bigr|
\end{eqnarray}
for some positive constant $C$. Noticing the elementary relation
%
\begin{eqnarray}
\bigl|\bigl(\Sigma^{1/2}\bigl(\mathcal{G}^{(1)}\bigr)^2
\Sigma^{1/2}\bigr)_{kk}\bigr| &\leq & \eta^{-1}\bigl|\Im \bigl(
\Sigma^{1/2}\bigl(\mathcal{G}^{(1)}\bigr)\Sigma^{1/2}
\bigr)_{kk}\bigr|
\nonumber
\\[-8pt]
\label{502.46}
\\[-8pt]
\nonumber
&\leq & \eta ^{-1}\bigl|\bigl(\Sigma^{1/2}
\mathcal{G}^{(1)}\Sigma^{1/2}\bigr)_{kk}\bigr|,
\end{eqnarray}
it thus suffices to bound the last term of (\ref{502.45}). In addition,
combining (\ref{502.46}) and (\ref{502.21}) we see
\[
\bigl|\bigl(\Sigma^{1/2}\bigl(\mathcal{G}^{(1)}\bigr)^2
\Sigma^{1/2}\bigr)_{kk}\bigr|\leq N^{1+C\varepsilon}
\]
holds with $\zeta$-high probability.
Finally, the estimate of the last term of (\ref{502.45}) can be
addressed as follows:
\begin{eqnarray*}
&&N^{-3}\sum_{k}\eta\bigl|\bigl(
\Sigma^{1/2}\mathcal{G}^{(1)}\Sigma ^{1/2}
\bigr)_{kk}\bigr|^2\bigl|\bigl(\Sigma^{1/2}\bigl(
\mathcal{G}^{(1)}\bigr)^2\Sigma ^{1/2}
\bigr)_{kk}\bigr|\\
&& \qquad \leq  N^{-7/3+C\varepsilon}\sum_{k}\bigl|
\bigl(\Sigma^{1/2}\mathcal {G}^{(1)}\Sigma^{1/2}
\bigr)_{kk}\bigr|
\leq N^{-7/3+C\varepsilon}\Trr \bigl(\Sigma^{1/2}\bigl|\mathcal{G}^{(1)}\bigr|
\Sigma ^{1/2}\bigr)\\
&&\qquad\leq C^\prime N^{-7/3+C\varepsilon}\Trr\bigl|
\mathcal{G}^{(1)}\bigr| =O\bigl(N^{-4/3+C\varepsilon}\bigr)
\end{eqnarray*}
holds with $\zeta$-high probability, where in the last step we used the
fact that $\Trr|\mathcal{G}^{(1)}|\leq N^{1+\varepsilon}$ with $\zeta
$-high probability for any fixed $\zeta>0$, which has been proved in
\cite{BPZ2013} (see Lemma~3.10 therein). Again, since $|(\Sigma
^{1/2}\mathcal{G}^{(1)}\Sigma^{1/2})_{ij}|$ are trivially bounded by
$O(\eta^{-1})$, the above bound also holds in expectation. Thus, we
complete the proof of (\ref{11.11.101}).
\end{pf*}
%
\section{Bounds on the entries of
\texorpdfstring{$\Sigma^{1/2}\mathcal{G}\Sigma^{1/2}$}
{Sigma{1/2}{G}Sigma{1/2}}}\label{sec5}

$\!$In this section, we prove Lem\-ma~\ref{lem.11.11.105}. Substantially
different strategies will be adopted for the proofs of (\ref{502.20})
and (\ref{502.21}). Thus we will perform them separately. Moreover,
since $\mathcal{G}^{(i)}$ and $\mathcal{G}$ are only different in
dimension (observing that they share the same population covariance
matrix $\Sigma$), we will harmlessly work on $\mathcal{G}$ for simplicity.
\begin{pf*}{Proof of (\ref{502.20}) (with $\mathcal{G}^{(1)}$ replaced
by $\mathcal{G}$)} Note when $\Sigma$ is diagonal, we can denote it as
$\Sigma=\operatorname{diag} (\sigma_1^2,\ldots, \sigma_M^2)$.
Let $\Sigma^{[j]}$ be the $(M-1)\times(M-1)$ minor of $\Sigma$,
obtained by deleting the $j$th column and row of $\Sigma$.
Moreover, we denote the $j$th row of $X$ by $\mathfrak{X}_j^*$ thus its
conjugate transpose by $\mathfrak{X}_j$, and denote by $X^{[j]}$ the
$(M-1)\times N$ submatrix obtained via deleting $\mathfrak{X}_j^*$ from
$X$. Correspondingly, we will use the notation
\begin{eqnarray*}
\mathcal{W}^{[j]} &=& \bigl(\Sigma^{[j]}\bigr)^{1/2}X^{[j]}
\bigl(X^{[j]}\bigr)^*\bigl(\Sigma ^{[j]}\bigr)^{1/2},\qquad
W^{[j]}=\bigl(X^{[j]}\bigr)^*\Sigma^{[j]}
X^{[j]},
\\
\mathcal{G}^{[j]} &=& \bigl(\mathcal{W}^{[j]}-z
\bigr)^{-1},\qquad G^{[j]}=\bigl(W^{[j]}-z
\bigr)^{-1}.
\end{eqnarray*}
Employing the Schur complement yields
\[
\mathcal{G}_{ii}=\frac{1}{-z-z\sigma_i^2\mathfrak{X}_i^*
((X^{[i]})^*\Sigma^{[i]} X^{[i]}-z )^{-1}\mathfrak{X}_i}=\frac
{1}{-z-z\sigma_1^2\mathfrak{X}_i^*G^{[i]}\mathfrak{X}_i}.
\]
Then by using Lemma~\ref{lem.k.40} and Theorem~\ref{thm.3.90} again, we
can actually get the following lemma, whose proof will be provided in
the supplementary material \cite{BPZ2014} in detail.

\begin{lem} \label{lem.2014050505050505} For any $\zeta>0$ given, there
exists some positive constant $C_\zeta$, such that for any $z\in
S_r(\tilde{c}, 5C_\zeta)$,
\[
\mathcal{G}_{ii}=1/\bigl(-z-z\sigma_i^2m_0(z)+o(1)
\bigr)
\]
holds with $\zeta$-high probability.
\end{lem}

Now we proceed to the proof of (\ref{502.20}). Ensured by (iii) of
Theorem~\ref{lem.2.2}, we deduce from Lemma~\ref{lem.2014050505050505}
that for $z\in S_r(\tilde{c}, 5C_\zeta)$,
%
\begin{equation}\label{x.x.x}
\bigl|\mathcal{G}_{ii}(z)\bigr|\leq C
\end{equation}
with $\zeta$-high probability for some positive constant $C$
independent of $z$. Therefore, we get the bound for $\mathcal{G}_{ii}$
when $\Sigma$ is diagonal. Actually we can strengthen (\ref{x.x.x}) to
the uniform bound as
%
\begin{equation}\label{1207.30}
\sup_{z\in S_r(\tilde{c}, 5C_\zeta)}\bigl|\mathcal{G}_{ii}(z)\bigr|=O(1)
\end{equation}
with $\zeta$-high probability. To see this, we can assign an
$\varepsilon$-net on the region $S_r(\tilde{c},5C_\zeta)$ with
$\varepsilon=N^{-100}$ (say). Then by the definition of $\zeta$-high
probability, we see that (\ref{x.x.x}) holds for all $z$ in this
$\varepsilon$-net uniformly with $\zeta$-high probability. Moreover,
note $|\mathcal{G}_{ii}'(z)|\leq N^2$ for $z\in S_r(\tilde{c},5C_\zeta
)$, thus by the Lipschitz continuity, we can extend the bound to the
whole region $S_r(\tilde{c}, 5C_\zeta)$ easily.

Now, we are ready to use (\ref{1207.30}) to derive the aforementioned
delocalization property for the eigenvectors of $\mathcal{W}$ in the
edge case. Then we use the delocalization result to bound $\mathcal
{G}_{ij}$ and $(\mathcal{G}^2)_{ij}$ in return. Denoting the unit
eigenvector of $\mathcal{W}$ corresponding to $\lambda_k(\mathcal{W})$ by
\[
\mathbf{u}_k=(u_{k1},\ldots,u_{kM})^T,
\]
we can formulate the following lemma.

\begin{lem}\label{lem.zzz}
When $\Sigma$ is diagonal, for $\lambda
_{k}(\mathcal{W})\in[\lambda_r-\tilde{c}/2,C_r]$, we have
\[
\max_i|u_{ki}|^2\leq
\varphi^{C_\zeta}N^{-1}
\]
with $\zeta$-high probability.
\end{lem}

\begin{pf}
By (\ref{1207.30}) and the spectral decomposition, we have
\[
\Im\mathcal{G}_{ii}(z)=\sum_{k=1}^M
\frac{\eta}{(\lambda_k(\mathcal
{W})-E)^2+\eta^2}|u_{ki}|^2=O(1),
\]
with $\zeta$-high probability.
Now we set $\eta=\varphi^{C_\zeta}N^{-1}$. In light of (\ref{1207.30}),
we can set $E=\lambda_k(\mathcal{W}$ if
$\lambda_k(\mathcal{W})\in[\lambda_r-\tilde{c}/2,C_r]$. Then with
$\zeta$-high probability,
\[
\frac{\eta}{(\lambda_k(\mathcal{W})-E)^2+\eta^2}|u_{ki}|^2= \varphi
^{-C_\zeta}N|u_{ki}|^2=O(1),
\]
which implies Lemma~\ref{lem.zzz} immediately. Thus, we complete the proof.
\end{pf}
Now relying on the above delocalization property, we proceed to prove
(\ref{502.20}).
Note that by the spectral decomposition, for $z$ satisfying the
assumption in Lemma~\ref{lem.11.11.105} and $\alpha=1,2$ we have
%
\begin{eqnarray}
\bigl|\bigl(\mathcal{G}^\alpha(z)\bigr)_{ij}\bigr| &\leq& \sum
_{k=1}^M\frac{1}{|\lambda
_k(\mathcal{W})-z|^\alpha}|u_{ki}||u_{kj}|
\nonumber\\
\label{201404271} &\leq & \sum_{k\dvtx \lambda_k\in
[\lambda_r-\tilde{c}/2,C_r]}\frac{1}{|\lambda_k-z|^\alpha
}|u_{ki}||u_{kj}|+O(1)
\\
\nonumber
&\leq& \varphi^{C_\zeta}\frac{1}N\sum
_{k=1}^M\frac{1}{|\lambda_k(\mathcal
{W})-z|^{\alpha}}+O(1)
\end{eqnarray}
with $\zeta$-high probability. When $\alpha=2$, we see that with $\zeta$-high probability,
\[
\frac{1}N\sum_{k=1}^M
\frac{1}{|\lambda_k(\mathcal{W})-z|^2}=\eta^{-1}\Im m_N(z)=\eta^{-1}
\biggl(\Im m_0(z)+O\biggl(\frac{\varphi^{C_\zeta}}{N\eta}\biggr)\biggr)
\]
according to (\ref{144201}) and (\ref{9.1.1}). From (ii) of Theorem~\ref
{lem.2.2} we have $\Im m_0(z)\sim\eta/\sqrt{\kappa+\eta}$. Noticing the
assumptions on $E$ and $\eta$ in Lemma~\ref{lem.11.11.105}, we
immediately get that
the second inequality in (\ref{502.20}) holds. Now, when $\alpha=1$, we
claim that for some sufficiently large constant $C_\zeta>0$,
%
\begin{equation}\label{1207.20}
\frac{1}N\sum_{k=1}^M
\frac{1}{|\lambda_k(\mathcal{W})-z|}=\frac{1}{N}\Trr \bigl|\mathcal{G}(z)\bigr|= O\bigl((\log
N)^{C_\zeta}\bigr)
\end{equation}
holds with $\zeta$-high probability.\vspace*{1.5pt} Such a bound has been established
in Lem\-ma~3.10 of \cite{BPZ2013} for $\frac{1}{N}\Trr |\mathcal
{G}^{(i)}(z)|$ by using the strong local MP type law. It is just the
same to check its validity for $\frac{1}{N}\Trr |\mathcal{G}(z)|$
[bearing in mind that for (\ref{502.20}) what we really need is the
bound for $\frac{1}{N}\Trr |\mathcal{G}^{(1)}(z)|$]. So we will not
reproduce the details here. Therefore, we complete the proof of (\ref{502.20}).
\end{pf*}

Now we start to tackle the much more complicated case, that is, (\ref
{502.21}) for general $\Sigma$.
\begin{pf*}{Proof of (\ref{502.21}) (with $\mathcal{G}^{(1)}$ replaced
by $\mathcal{G}$)} For simplicity, we will also work with $\mathcal{G}$
instead of $\mathcal{G}^{(1)}$.
Note that
%
\begin{equation}\label{502.120}
\quad\Sigma^{1/2}\mathcal{G}\Sigma^{1/2}=\Sigma^{1/2}
\bigl(\Sigma ^{1/2}XX^*\Sigma^{1/2}-zI \bigr)^{-1}
\Sigma^{1/2}= \bigl(XX^*-z\Sigma ^{-1} \bigr)^{-1}.
\end{equation}
For convenience, we use the notation $\Phi:=\Sigma^{-1}$ and recall
$\Delta=\Delta(z):=\Sigma^{1/2}\mathcal{G}(z)\Sigma^{1/2}$
defined in \hyperref[sec1]{Introduction}.
Thus, we have $\Delta_{kk}:=\Delta_{kk}(z)=\break [(XX^*-z\Phi)^{-1}]_{kk}$.
An elementary observation from the spectral decomposition is
%
\begin{equation}\label{2014050311}
\Delta_{kk}=\sum_{i=1}^M
\frac{1}{\lambda_i(\mathcal{W})-z}\bigl(\Sigma ^{1/2}\mathbf{u}_i
\mathbf{u}_i^*\Sigma^{1/2}\bigr)_{kk}.
\end{equation}
We will only provide the estimate for $\Delta_{11}$ in the sequel,
since the others can be handled analogously. The following lemma lies
at the core of our subsequent discussion.

\begin{lem} \label{lem.502.70}
Let $z_0:=E_0+\mathbf{i}\eta_0$ satisfy
$E_0\in[\lambda_r-\tilde{c}, \lambda_r+N^{-2/3+\varepsilon}]$ and $\eta
_0:=N^{-2/3+A_0\varepsilon}$
for some positive constant $A_0>1$ independent of $\varepsilon$. Under
Condition~\ref{con.1.1}, for any given constant $\zeta>0$ we have
%
\begin{equation}\label{502.100}
\sup_{E_0\in[\lambda_r-\tilde{c}, \lambda_r+N^{-2/3+\varepsilon
}]}\bigl|\Delta_{11}(z_0)\bigr|\leq C
\eta_0^{-1/2}
\end{equation}
with $\zeta$-high probability for some positive constant $C$.
\end{lem}

We postpone the proof of Lemma~\ref{lem.502.70} to the end of this
section and proceed to prove (\ref{502.21}) by assuming Lemma~\ref
{lem.502.70}. By (\ref{502.100}) and the spectral decomposition we have
%
\begin{equation}\label{502.101}
C\eta_0^{-1/2}\geq\Im\Delta_{11}(z_0)=
\sum_i\frac{\eta_0}{(\lambda
_i(\mathcal{W})-E_0)^2+\eta_0^2}\bigl(
\Sigma^{1/2}\mathbf{u}_i\mathbf {u}_i^*
\Sigma^{1/2}\bigr)_{11}
\end{equation}
with $\zeta$-high probability. We set in (\ref{502.101}) that
$E_0=\lambda_i(\mathcal{W})$ for some $\lambda_i(\mathcal{W})\in
[\lambda_r-\tilde{c},\lambda_r+N^{-2/3+\varepsilon}]$. Immediately,
(\ref{502.101}) implies that
%
\begin{equation}\label{502.102}
\bigl(\Sigma^{1/2}\mathbf{u}_i\mathbf{u}_i^*
\Sigma^{1/2}\bigr)_{11}\leq C\eta _0^{1/2}
\leq N^{-1/3+A_0\varepsilon}
\end{equation}
holds with $\zeta$-high probability.
(\ref{502.102}) together with (\ref{502.100}) can then be employed to
bound $\Delta_{11}(z)$, for all $z$ satisfying the assumption of Lemma~\ref{lem.11.11.105}. We perform it as follows. At first, according to
(i) of Theorem~\ref{thm.1206.7}, we can assume that $\lambda_1(\mathcal
{W})\leq\lambda_r+N^{-2/3+\varepsilon}$. Now, for $z=E+\mathbf{i}\eta
$, we choose $E_0=E$ thus $z_0=E+\mathbf{i}\eta_0$. Again, by the
spectral decomposition, we see that
%
\begin{eqnarray}
\bigl|\Delta_{11}(z)-\Delta_{11}(z_0)\bigr|
&=& \Biggl|\sum
_{i=1}^M\biggl(\frac{1}{\lambda
_i(\mathcal{W})-z}-
\frac{1}{\lambda_i(\mathcal{W})-z_0}\biggr) \bigl(\Sigma ^{1/2}\mathbf{u}_i
\mathbf{u}_i^*\Sigma^{1/2}\bigr)_{11}\Biggr|
\nonumber
\\
&\leq & (\eta_0-\eta)\sum_{i=1}^M
\frac{(\Sigma^{1/2}\mathbf{u}_i\mathbf
{u}_i^*\Sigma^{1/2})_{11}}{|(\lambda_i(\mathcal{W})-z)(\lambda
_i(\mathcal{W})-z_0)|}\nonumber\\
\label{2014050312}&\leq & (\eta_0-\eta)\sum
_{i=1}^M \frac{(\Sigma
^{1/2}\mathbf{u}_i\mathbf{u}_i^*\Sigma^{1/2})_{11}}{|\lambda_i(\mathcal
{W})-z|^2}
\\
&=&  (\eta_0-\eta)\eta^{-1}\sum
_{i=1}^M \Im\frac{(\Sigma^{1/2}\mathbf
{u}_i\mathbf{u}_i^*\Sigma^{1/2})_{11}}{\lambda_i(\mathcal{W})-z}\nonumber\\
&\leq &  N^{2A_0\varepsilon}
\sum_{i=1}^M \Im\frac{(\Sigma^{1/2}\mathbf
{u}_i\mathbf{u}_i^*\Sigma^{1/2})_{11}}{\lambda_i(\mathcal{W})-z}
\nonumber
\end{eqnarray}
with $\zeta$-high probability. Now we split the index collection $\{
1,\ldots,M\}$ into two parts as
\[
I_1:=\bigl\{i\dvtx  \lambda_i(\mathcal{W})\in\bigl[
\lambda_r-\tilde{c},\lambda _r+N^{-2/3+\varepsilon}\bigr]\bigr\}, \qquad I_2:=\bigl\{i\dvtx \lambda_i(\mathcal{W})< \lambda
_r-\tilde{c}\bigr\}.
\]
Combining (\ref{502.102}), (\ref{2014050312}) and the assumption on $z$ yields
\begin{eqnarray*}
&& \bigl|\Delta_{11}(z)-\Delta_{11}(z_0)\bigr|\\
 &&\qquad\leq
N^{-1/3+3A_0\varepsilon}\sum_{i\in I_1}\Im\frac{1}{\lambda_i(\mathcal{W})-z}+CN^{2A_0\varepsilon
}
\eta\sum_{i\in I_2}\bigl(\Sigma^{1/2}
\mathbf{u}_i\mathbf{u}_i^*\Sigma ^{1/2}
\bigr)_{11}
\\
&& \qquad\leq   N^{2/3+3A_0\varepsilon}\Im m_N(z)+CN^{2A_0\varepsilon}\eta
\Sigma_{11}\\
&&\qquad \leq N^{2/3+4A_0\varepsilon}\biggl(\Im m_0(z)+
\frac{1}{N\eta
}\biggr)+CN^{-2/3+2A_0\varepsilon}\Sigma_{11}
\\
&&\qquad\leq   N^{1/3+5A_0\varepsilon}
\end{eqnarray*}
with $\zeta$-high probability. Here in the last two inequalities we
used (\ref{9.1.1}) and (ii) of Theorem~\ref{lem.2.2}, along with the
fact that $\Sigma_{11}$ is bounded. Hence, we have
%
\begin{equation}\label{201404272}
\bigl|\Delta_{11}(z)\bigr|\leq\bigl|\Delta_{11}(z_0)\bigr|+N^{1/3+5A_0\varepsilon}
\leq N^{1/3+6A_0\varepsilon}
\end{equation}
with $\zeta$-high probability.
Thus (\ref{502.21}) follows if we replace $\mathcal{G}$ by $\mathcal{G}^{(1)}$.
\end{pf*}

The remaining part of this section will be devoted to the proof of
Lem\-ma~\ref{lem.502.70}.
\begin{pf*}{Proof of Lemma~\ref{lem.502.70}}
At first, analogous to
the derivation of (\ref{1207.30}) via (\ref{x.x.x}), the verification
of (\ref{502.100}) can be reduced to providing the desired bound on
$|\Delta_{11}(z_0)|$ for any single $z_0$ with $E_0\in[\lambda_r-\tilde
{c}, \lambda_r+N^{-2/3+\varepsilon}]$. Hence, in the sequel, we will
just fix $E_0$. The extension to the uniform bound via Lipschitz
continuity and $\varepsilon$-net argument is just routine. We recall
the notation $ \mathfrak{X}_j$ and $X^{[j]}$ in the proof of (\ref{502.20}).
For simplicity, we further write
\[
\Sigma^{-1}=\Phi:=\pmatrix{
\phi_{11} &\Phi_1^*
\vspace*{3pt}\cr
\Phi_1 &\Phi^{[1]}},
\]
where $\phi_{11}$ is the $(1,1)$th entry of $\Phi$ and $\Phi_1$ is its
first column with $\phi_{11}$ removed. As the inverse of $\Sigma$, we
know that $\Phi$ is also positive-definite and its eigenvalues are
bounded both from below and above, in light of Condition \ref{con.1.1}.
Consequently, its entries are also bounded, so is $\Vert\Phi_1\Vert$. Now by
using Schur complement to (\ref{502.120}) we can
deduce that
\begin{eqnarray*}
&& \Delta_{11}(z_0)\\
&&\qquad={1}/\bigl(\mathfrak{X}_1^*\mathfrak{X}_1-z_0\phi
_{11}\\
&&\hspace*{3pt}\qquad\qquad{}- \bigl(\mathfrak{X}_1^*\bigl(X^{[1]}\bigr)^*-z_0\Phi_1^* \bigr)
\bigl(X^{[1]}\bigl(X^{[1]}\bigr)^*-z_0\Phi^{[1]} \bigr)^{-1} \bigl(X^{[1]}\mathfrak
{X}_1-z_0\Phi_1 \bigr)\bigr)
\\
&&\qquad:=\frac{1}{D_1+D_2+D_3},
\end{eqnarray*}
where $D_i:=D_i(z_0),i=1,2,3$, whose explicit formulas are as follows,
\begin{eqnarray*}
D_1 &:=& \mathfrak{X}_1^*\mathfrak{X}_1-z_0
\phi_{11}-\mathfrak {X}_1^*\bigl(X^{[1]}
\bigr)^* \bigl(X^{[1]}\bigl(X^{[1]}\bigr)^*-z_0
\Phi^{[1]} \bigr)^{-1}X^{[1]}\mathfrak{X}_1,
\\
D_2 &:=& -z_0^2\Phi_1^*
\bigl(X^{[1]}\bigl(X^{[1]}\bigr)^*-z_0
\Phi^{[1]} \bigr)^{-1}\Phi_1,
\\
D_3 &:=& z_0\mathfrak{X}_1^*
\bigl(X^{[1]}\bigr)^* \bigl(X^{[1]}\bigl(X^{[1]}
\bigr)^*-z_0\Phi ^{[1]} \bigr)^{-1}
\Phi_1\\
&&{}+z_0\Phi_1^* \bigl(X^{[1]}
\bigl(X^{[1]}\bigr)^*-z_0\Phi ^{[1]}
\bigr)^{-1}X^{[1]}\mathfrak{X}_1.
\end{eqnarray*}
Our starting point is the following elementary inequality:
%
\begin{equation}\label{1207.40}
\hspace*{9pt}\bigl|\Delta_{11}(z_0)\bigr|\leq\min\bigl\{\bigl(\bigl|
\Im(D_1+D_2+D_3)\bigr|\bigr)^{-1},\bigl|\Re
(D_1+D_2+D_3)\bigr|^{-1}\bigr\}.
\end{equation}
Observe that if $|\Re(D_1+D_2+D_3)|> N^{1/6}$, the bound for $|\Delta
_{11}(z_0)|$ in (\ref{502.100}) automatically holds.
Hence, it suffices to show that with $\zeta$-high probability,
%
\begin{equation}\label{502.300}
\bigl|\Im(D_1+D_2+D_3)\bigr|\geq C
\eta_0^{1/2}
\end{equation}
when
%
\begin{equation}\label{502.301}
\bigl|\Re(D_1+D_2+D_3)\bigr|\leq N^{1/6}
\end{equation}
for some positive constant $C$. In order to verify (\ref{502.300})
under assumption (\ref{502.301}), a careful analysis on the real and
imaginary parts of $D_1,D_2,D_3$ is required. We perform it as follows.
We start from the following reduction on $D_1$,
\begin{eqnarray*}
D_1&=&\mathfrak{X}_1^*\mathfrak{X}_1-z_0
\phi_{11}\\
&&{}-\mathfrak {X}_1^*\bigl(X^{[1]}
\bigr)^*\bigl(\Phi^{[1]}\bigr)^{-{1}/2} \bigl(\bigl(
\Phi^{[1]}\bigr)^{-{1}/2}X^{[1]}\bigl(X^{[1]}
\bigr)^*\bigl(\Phi^{[1]}\bigr)^{-{1}/2}-z_0
\bigr)^{-1}\\
&&\hspace*{11pt}{}\times\bigl(\Phi ^{[1]}\bigr)^{-{1}/2}X^{[1]}
\mathfrak{X}_1
\\
&=& \mathfrak{X}_1^*\mathfrak{X}_1-z_0
\phi_{11}-\mathfrak {X}_1^*\bigl(X^{[1]}
\bigr)^*\bigl(\Phi^{[1]}\bigr)^{-1}X^{[1]} \bigl(
\bigl(X^{[1]}\bigr)^{*}\bigl(\Phi ^{[1]}
\bigr)^{-1}X^{[1]}-z_0 \bigr)^{-1}
\mathfrak{X}_1
\\
&=&-z_0\phi_{11}-z_0\mathfrak{X}_1^*
\bigl(\bigl(X^{[1]}\bigr)^{*}\bigl(\Phi ^{[1]}
\bigr)^{-1}X^{[1]}-z_0 \bigr)^{-1}
\mathfrak{X}_1.
\end{eqnarray*}
In the second equality above, we have used the elementary fact that for
any $m\times n$ matrix $A$
\[
A\bigl(A^*A-z_0I_n\bigr)^{-1}A^*=AA^*
\bigl(AA^*-z_0I\bigr)^{-1},
\]
which can be checked by the singular decomposition easily.
To abbreviate, we use the notation
\[
\widetilde{G}^{[1]}(z_0):= \bigl(\bigl(X^{[1]}
\bigr)^{*}\bigl(\Phi ^{[1]}\bigr)^{-1}X^{[1]}-z_0
\bigr)^{-1}.
\]
Adopting Lemma~\ref{lem.k.40} again, we obtain
%
\begin{equation}\label{201405110101010}
D_1=-z_0\phi_{11}-z_0
\frac{1}N\Trr\widetilde{G}^{[1]}(z_0)+O \biggl(
\frac
{\varphi^{C_\zeta}}{N}\bigl\Vert\widetilde{G}^{[1]}(z_0)\bigr\Vert_{\mathrm{HS}}
\biggr)
\end{equation}
with $\zeta$-high probability.
Now, we need the following lemma whose proof will be also stated in the
supplementary material \cite{BPZ2014}.

\begin{lem} \label{lem.1207.50}
Under the above notation, we can show that
%
\begin{equation}\label{1208.2}
\frac{1}N \Trr \widetilde{G}^{[1]}(z_0)=m_N(z_0)+O
\biggl(\frac{1}{N\eta_0}\biggr).
\end{equation}
\end{lem}

Denoting $\kappa_0:=|\lambda_r-E_0|$, we deduce from Lemma~\ref
{lem.1207.50} that
%
\begin{equation}\label{502.405}
\biggl|\frac{1}N \Trr \widetilde{G}^{[1]}(z_0)\biggr|=O(1),\qquad
\frac{1}{N}\Im \Trr \tilde{G}^{[1]}(z_0)\sim\sqrt{
\kappa_0+\eta_0}
\end{equation}
hold with $\zeta$-high probability, by combining (\ref{9.1.1}) and
(i)--(ii) of Theorem~\ref{lem.2.2}.
By~(\ref{144201}), we have
$
\Vert\widetilde{G}^{[1]}(z_0)\Vert_{\mathrm{HS}}=\sqrt{\Im \Trr \widetilde
{G}^{[1]}(z_0)/\eta_0}$,
which together with (\ref{201405110101010}) and (\ref{502.405}) implies that
%
\begin{equation}\label{502.550}
\bigl|\Re D_1(z_0)\bigr|\leq\bigl|D_1(z_0)\bigr|=O(1)
\end{equation}
and
%
\begin{equation}\label{1207.91}
\hspace*{15pt}\Im D_1(z_0)=-E_0\frac{1}{N}\Im
\Trr \widetilde{G}^{[1]}(z_0)+O \biggl(
\varphi^{C_\zeta}\sqrt{\frac{\Im \Trr \widetilde{G}^{[1]}(z_0)}{N^2\eta
_0}} \biggr)+O(\eta_0)
\end{equation}
with $\zeta$-high probability. Here, we also used the fact that $|\phi
_{11}|$ is bounded. Then by (\ref{1208.2}), (\ref{1207.91}) and (\ref
{9.1.1}) we have
%
\begin{equation}\label{502.801}
\Im D_1=-E_0\Im m_0(z_0)+O
\bigl(N^{-1/3-C\varepsilon}\bigr)
\end{equation}
with $\zeta$-high probability.

We proceed to the analysis toward $D_2$ and $D_3$. For $D_2$, by
definition we have
%
\begin{eqnarray}
\Re D_2(z_0)&=&-\bigl(E_0^2-
\eta_0^2\bigr)\Re\Phi_1^*
\bigl(X^{[1]}\bigl(X^{[1]}\bigr)^*-z_0
\Phi^{[1]} \bigr)^{-1}\Phi_1
\nonumber
\\[-8pt]
\label{502.551}
\\[-8pt]
\nonumber
&&{}+2E_0\eta_0\Im\Phi_1^*
\bigl(X^{[1]}\bigl(X^{[1]}\bigr)^*-z_0
\Phi^{[1]} \bigr)^{-1}\Phi_1,
\\
\Im D_2(z_0)&=&-\bigl(E_0^2-
\eta_0^2\bigr)\Im\Phi_1^*
\bigl(X^{[1]}\bigl(X^{[1]}\bigr)^*-z_0
\Phi^{[1]} \bigr)^{-1}\Phi_1
\nonumber
\\[-8pt]
\label{502.600}
\\[-8pt]
\nonumber
&& {}-2E_0\eta_0\Re\Phi_1^*
\bigl(X^{[1]}\bigl(X^{[1]}\bigr)^*-z_0
\Phi^{[1]} \bigr)^{-1}\Phi_1.
\end{eqnarray}
Now, for $D_3$, we have the following lemma whose proof will be
presented in the supplementary material \cite{BPZ2014}.

\begin{lem} \label{lem.20140505050505050505050}
Assume that $z_0$ satisfies the assumption in Lemma~\ref{lem.502.70}.
For any $\zeta>0$, there exists some constant $C_\zeta$ such that
%
\begin{equation}\label{502.552}
|D_3|\leq\frac{\varphi^{C_\zeta}}{\sqrt{N}}\sqrt{\eta_0^{-1}
\Im\Phi _1^* \bigl(X^{[1]}\bigl(X^{[1]}
\bigr)^*-z_0\Phi^{[1]} \bigr)^{-1}
\Phi_1}
\end{equation}
holds with $\zeta$-high probability.
\end{lem}

Now we invoke the crude bound
%
\begin{equation}\label{1207.90}
\hspace*{10pt}\Im\Phi_1^* \bigl(X^{[1]}\bigl(X^{[1]}
\bigr)^*-z_0\Phi^{[1]} \bigr)^{-1}\Phi
_1\leq C\eta_0^{-1}\bigl\Vert\bigl(
\Phi^{[1]}\bigr)^{-1/2}\Phi_1\bigr\Vert \leq
C_1\eta_0^{-1}
\end{equation}
with some positive constants $C$ and $C_1$, which trivially implies that
%
\begin{equation}\label{2014050570000001}
|D_3|\leq\frac{C\varphi^{C_\zeta}}{\sqrt{N}}\eta _0^{-1}=O
\bigl(N^{1/6-C_2\varepsilon}\bigr)
\end{equation}
with $\zeta$-high probability for some positive constant $C_2$. In
addition, plugging (\ref{1207.90}) into (\ref{502.551}) yields that
%
\begin{equation}\label{1207.100}
\hspace*{7pt}\Re D_2(z_0)= -\bigl(E_0^2-
\eta_0^2\bigr)\Re\Phi_1^*
\bigl(X^{[1]}\bigl(X^{[1]}\bigr)^*-z_0
\Phi^{[1]} \bigr)^{-1}\Phi_1+O(1)
\end{equation}
with $\zeta$-high probability.
Now,\vspace*{1pt} we are ready to provide a bound for $\Re\Phi_1^*
(X^{[1]}\times\break (X^{[1]})^*-z_0\Phi^{[1]} )^{-1}\Phi_1$ which is needed to
estimate $\Im D_2$ according to~(\ref{502.600}). Combining (\ref
{502.550}), (\ref{2014050570000001}) and (\ref{1207.100}), we can see that
\[
\bigl|\Re(D_1+D_2+D_3)\bigr|= \bigl|\bigl(E_0^2-
\eta_0^2\bigr)\Re\Phi_1^*
\bigl(X^{[1]}\bigl(X^{[1]}\bigr)^*-z_0
\Phi^{[1]} \bigr)^{-1}\Phi _1\bigr|+O
\bigl(N^{1/6-C_2\varepsilon}\bigr)
\]
with $\zeta$-high probability.
Now, invoking assumption (\ref{502.301}), we obtain
%
\begin{equation}\label{502.601}
\bigl|\Re\Phi_1^* \bigl(X^{[1]}\bigl(X^{[1]}
\bigr)^*-z_0\Phi^{[1]} \bigr)^{-1}\Phi
_1\bigr|=O\bigl(N^{1/6}\bigr)
\end{equation}
with $\zeta$-high probability.
Inserting (\ref{502.601}) into (\ref{502.600}) we have
%
\begin{equation}\label{502.800}
\hspace*{2pt}\qquad\Im D_2=-\bigl(E_0^2-\eta_0^2
\bigr)\Im\Phi_1^* \bigl(X^{[1]}\bigl(X^{[1]}
\bigr)^*-z_0\Phi ^{[1]} \bigr)^{-1}
\Phi_1+O\bigl(N^{-1/2+C\varepsilon}\bigr).
\end{equation}
For convenience, we set $t_0=\Im\Phi_1^* (X^{[1]}(X^{[1]})^*-z_0\Phi
^{[1]} )^{-1}\Phi_1$. Putting (\ref{502.552}), (\ref{502.800}) and
(\ref{502.801}) together, we get
\begin{eqnarray*}
&& \Im(D_1+D_2+D_3)\\
&& \qquad= -E_0\Im
m_0(z_0)-\bigl(E_0^2-
\eta_0^2\bigr)t_0+O \biggl(
\frac
{\varphi^{C_\zeta}}{\sqrt{N\eta_0}}t_0^{1/2} \biggr)+O\bigl(N^{-1/3-C\varepsilon}
\bigr)
\end{eqnarray*}
with $\zeta$-high probability. Now observe that $E_0\Im m_0(z_0)$ and
$(E_0^2-\eta_0^2)t_0$ are both positive. Moreover, by (ii) of Theorem~\ref{lem.2.2} we see that
%
\begin{equation}\label{1207.96}
\Im m_0(z_0)\sim\sqrt{\kappa_0+
\eta_0}.
\end{equation}
Now we split the discussion into two cases according to whether
%
\begin{equation}\label{1207.92}
t_0\gg\frac{\varphi^{C_\zeta}}{\sqrt{N\eta_0}}t_0^{1/2},
\end{equation}
holds. If (\ref{1207.92}) is valid, then we deduce from (\ref{1207.96})
that (\ref{502.300}) holds. If (\ref{1207.92}) fails, we claim that one
must have
%
\begin{equation}\label{1207.95}
\frac{\varphi^{C_\zeta}}{\sqrt{N\eta_0}}t_0^{1/2}\ll\sqrt{\kappa_0+
\eta _0}.
\end{equation}
Since if (\ref{1207.92}) does not hold, there exists some positive
constant $C$ such that $t_0\leq C\varphi^{2C_\zeta}/N\eta_0$,
which implies (\ref{1207.95}) immediately by our choice of $\eta_0$.
Now (\ref{1207.96}) and (\ref{1207.95}) imply (\ref{502.300}) again.
Then by (\ref{1207.40}), we complete the proof.
\end{pf*}

\section*{Acknowledgements}
The authors are very grateful to Laszlo Erd\"{o}s, Horng-Tzer Yau and
Jun Yin for useful comments on a former version of this paper which
helped to improve the organization significantly. We also thank the
anonymous referees for careful reading and valuable comments,
especially the suggestions on statistical signal processing.

\begin{supplement}[id=suppA]
\sname{Supplement}
\stitle{Proofs of some lemmas}
\slink[doi]{10.1214/14-AOS1281SUPP} 
\sdatatype{.pdf}
\sfilename{aos1281\_supp.pdf}
\sdescription{In the supplementary material \cite{BPZ2014}, we will
provide the proofs of Lemmas \ref{lem.2014050501}, \ref{lem.y.4.6}, \ref
{lem.20140505123}, \ref{lem.2014050505050505}, \ref{lem.1207.50} and
\ref{lem.20140505050505050505050}.}
\end{supplement}




\printaddresses
\end{document}